\let\ORIlabel\label
\let\ORIrefstepcounter\refstepcounter
\AddToHook{package/hyperref/before}{
   \let\label\ORIlabel 
   \let\refstepcounter\ORIrefstepcounter}

\documentclass{siamart220329}

\usepackage{tgadventor}
\usepackage{url}
\usepackage{amsmath,amssymb}
\usepackage{graphicx,subcaption}
\usepackage{multirow,float,tikz,pgfmath}
\usetikzlibrary{patterns.meta,arrows,snakes,backgrounds,patterns}
\usepackage{algorithm}
\usepackage{algpseudocode}

\newsiamremark{remark}{Remark}
\crefname{remark}{remark}{remarks}
\Crefname{remark}{Remark}{Remarks}

\newcommand{\R}[1]{\mathbb{R}^{#1}}
\newcommand{\HT}{$\mathcal{H}^{2}$~}
\newcommand{\HTMG}{$\mathcal{H}^{2}$-MG~}
\newcommand{\OB}{\mathcal{O}}

\definecolor{darkred}{RGB}{173,0,43}
\definecolor{dgreen}{RGB}{0, 115, 85}
\definecolor{lgreen}{RGB}{217, 231, 225}
\definecolor{lblue}{RGB}{0, 115, 85} 
\definecolor{newblue}{RGB}{0,0,0}

\title{$\mathcal{H}^{2}$-MG: A multigrid method for hierarchical rank structured matrices}
\author{
Daria Sushnikova\thanks{King Abdullah University of Science and Technology (KAUST), Thuwal, Saudi Arabia.}
\and George Turkiyyah\footnotemark[1]
\and Edmond Chow\thanks{School of Computational Science and Engineering, Georgia Institute of Technology,  Atlanta, GA.  The work of this author was supported by the U.S. National Science Foundation, award OAC-2003683.}
\and David Keyes\footnotemark[1]
}
\date{\today}

\begin{document}
\usetikzlibrary {arrows.meta}
\maketitle

\begin{abstract}
This paper presents a new fast iterative solver for large systems involving kernel matrices. Advantageous aspects of \HT matrix approximations and the multigrid method are hybridized to create the \HTMG algorithm. This combination provides the time and memory efficiency of \HT operator representation along with the rapid convergence of a multilevel method. 
We describe how  \HTMG works, show its linear complexity, and demonstrate its effectiveness on two standard kernels {\color{black}and on a single-layer potential boundary element discretization with complex geometry}. 
The current zoo of \HT solvers, which includes a wide variety of iterative and direct solvers, so far lacks a method that exploits multiple levels of resolution, commonly referred to in the iterative methods literature as ``multigrid'' from its origins in a hierarchy of grids used to discretize differential equations.  This makes \HTMG a valuable addition to the collection of \HT solvers.  The algorithm has potential for advancing various fields that require the solution of large, dense, symmetric positive definite matrices.
\end{abstract} 

\begin{keywords}
$\mathcal{H}^{2}$-matrix, multigrid methods,
kernel matrices, 
rank-structured matrices, iterative solvers, linear complexity 
\end{keywords}

\begin{MSCcodes}
65F10, 65N55, 65F30, 65F55
\end{MSCcodes}

\section{Introduction}

This paper tackles the challenge of solving linear systems with large, dense kernel matrices. Such systems arise in a wide range of applications, including computational statistics~\cite{bane-app_st-2005,will-app_st-2006}, machine learning~\cite{scho-app_ml-2002,rahi-app_ml-2007}, and computational physics~\cite{barn-app_ph-1986,doum-app_ph-2024}. Solving these systems is particularly challenging due to their quadratic and cubic complexity in terms of memory and runtime, respectively. Over the past decades, significant progress has been made to address this issue through rank-structured matrix approximations. 

Rank-structured methods, and particularly $\mathcal{H}^{2}$ matrices \cite{hack-h2-2000,borm-h2-2010}, typically provide a time- and memory-efficient matrix-vector product \cite{boukaram19a}, which naturally leads to solving systems of equations using iterative solvers. However, iterative solvers have their disadvantages, as their efficiency depends on the number of iterations and thus on the matrix conditioning. Direct solvers have their own challenges, being extremely complex for rank-structured matrices and involving a large constant overhead. To fill this gap, we consider the multigrid method \cite{fedo-mg-1964,bran-mg-1977,hack-mg-2013}, which is typically used for sparse matrices, and adapt it to \HT matrices. Multigrid methods exhibit excellent convergence properties and can significantly benefit from fast \HT matrix-vector products. In this paper, we introduce a new algorithm, $\mathcal{H}^{2}$-multigrid ($\mathcal{H}^{2}$-MG), which leverages the hierarchical structure of \HT matrices to create a multigrid method that operates across different levels of the \HT matrix. 

The key contributions of this work include:
\begin{itemize}
    \item Developing the \HTMG algorithm, a hybrid solver combining multilevel resolution with a hierarchical matrix representation.
    \item Demonstrating the linear complexity of \HTMG in both time and memory.
    \item Validating the effectiveness of \HTMG on problems from two standard kernel functions and a boundary element method, and comparing its performance with existing approaches.
\end{itemize}

The proposed method merges the time and memory efficiency of \HT matrices with the fast convergence properties of multigrid. This paper not only expands the repertoire of $\mathcal{H}^{2}$ solvers but also addresses a critical gap in the literature by providing a multigrid-inspired approach to hierarchical matrix methods. The simplicity and scalability of the $\mathcal{H}^{2}$-MG algorithm make it a valuable addition to the field, with potential applications across diverse domains.

\section{Related work}

A general $N$ by $N$ matrix requires $\OB(N^2)$ operations to compute a matrix-vector product and $\OB(N^2)$ memory for storage. However, significant progress has been made to reduce this computational cost with controllable loss of accuracy. Block low-rank matrix representations, such as the mosaic skeleton~\cite{tyrt-ms-1996,tyrt-ms-2000}, $\mathcal{H}$-matrices~\cite{khor-h-2000,hackb-h-2015,le-h-2006}, and HODLR matrices~\cite{ambi-hodlr-2013}, reduce the number of operations for matrix-vector products and storage to $\OB(N\log N)$ by exploiting low-rank approximations of certain blocks of the matrix. Further, nested-basis representations like HSS matrices~\cite{xia-hss-2010,chan-hss-2006,gill-hss-2012}, \HT matrices~\cite{hack-h2-2000,borm-h2-2010,mikh-mcbh-2016}, and the Fast Multipole Method (FMM)~\cite{grro-fmm-1987,ying-fmm-2004,gree-fmm-2021} were proposed. These methods enable matrix-vector multiplication with high accuracy in $\OB(N)$ operations for many matrices arising in physically causal models where interactions decay smoothly with distance \cite{boukaram19a,boukaram19b}. HSS and HODLR methods are particularly efficient for 1D problems, while \HT and FMM extend their efficiency to 2D and 3D problems thanks to their strong admissibility property.

This accelerated matrix-vector multiplication becomes the basis for solving linear systems using iterative techniques such as GMRES~\cite{saad-gmres-1986}, CG~\cite{hest-cg-1952,shew-cg-1994}, BiCGstab~\cite{van-bisgstab-1992}, and other iterative solvers. Iterative methods, while versatile, have a drawback: their convergence rate depends on the conditioning or eigenvalue clustering of the matrix. Thus, preconditioning for iterative methods, where the matrix is in \HT-matrix format, is often a necessity~\cite{xing-h2_precon-2021,zhao-afn-2024}. In contrast, fast direct solvers guarantee a predefined number of operations for solving the system. Although the complexity for general dense matrices remains $\OB(N^3)$, leveraging the hierarchical matrix format has led to breakthroughs in direct solver efficiency. Researchers in~\cite{hack-hlu-2000} introduced an $\OB(N\log N)$ direct solver algorithm, albeit with a substantial constant. Subsequent works~\cite{ambi-ifmm-2014,mind-rs-2017,sush-fmmlu-2023,ma-h2dir-2019,ma-h2dir-2024,Bouk-h2dir-2025,yesy-h2dir-2023} utilizing \HT (FMM) format have achieved $\OB(N)$ direct solver algorithms with more favorable constants.  

In this paper, we expand the zoo of \HT solvers by introducing a novel iterative algorithm---the \HTMG solver. This solver is rooted in the standard multigrid method but tailored for \HT structures. 

 Multigrid (MG) {\color{newblue}\cite{bran-mg-1977,hack-mg-2013}}  methods, as well as rank structured methods, play a crucial role in solving large linear systems, particularly systems with sparse matrices arising from discretized partial differential equations. 
Multigrid has its origins in~\cite{fedo-mg-1964} and was further developed in~\cite{bran-mg-1977}. Multigrid has become a cornerstone of numerical techniques for solving large sparse linear systems. Its strategy lies in efficiently reducing errors at multiple scales by leveraging a hierarchy of coarser grids. Over time, a variety of fruitful generalizations have emerged, such as algebraic multigrid (AMG)~\cite{bran-amg-1984} approaches that extend its applicability to unstructured grids and irregular geometries~\cite{falg-mg4-2002,de-mg5-2010,gref-mg11-2023} and fast multipole preconditioners for sparse matrices~\cite{ibei-mg12-2018}. Multigrid methods have also been turned into \HT matrices in order to provide fast methods for evaluating integral operators \cite{borm04_h2}. 

Our approach of applying multigrid to the \HT matrix diversifies the landscape of \HT solvers, offering a promising alternative that is easier to implement and parallelize efficiently compared to direct solvers. Demonstration of parallel scaling, however, is beyond the scope of this initial description.

\section{Algorithm}
\subsection{\HT matrix}
In this section, we briefly review the fundamental concept of an \HT matrix, a hierarchical block low-rank matrix structure with a nested basis property. Matrices well approximated in \HT form typically come from the discretization of boundary integral equations and several other problems with approximately separable kernels. 
The hierarchical nature of \HT matrices is the main inspiration for the \HTMG algorithm. For comprehensive and formal \HT definitions, see~\cite{grro-fmm-1987,hack-h2-2000,borm-h2-2010}.

Consider the linear system 
$$Ax = b,$$ where $A\in \R{N\times N}$, is dense and $x,b\in \R{N}$.
Let the rows and columns of matrix $A$ be partitioned into $M$ blocks. The size of $i$-th block is $B_i$, $i\in 1,\dots.,M$. 
Each block $A_{ij}$, $i,j\in 1,\dots.,M$, of matrix $A$ has either full rank, denoted by
$$A_{ij} = \widehat{D}_{ij}, \quad \widehat{D}_{ij} \in \R{B_i, B_j},$$
and is called a ``close'' block, or has a low rank, is called ``far'',
and possesses the following property: 
\begin{equation}
    A_{ij} \approx \widehat{F}_{ij} = \widehat{U}_i \widehat{S}_{ij} \widehat{V}_j,
    \label{eq:f}
\end{equation}
\noindent
with $\widehat{U}_i\in\R{B_i \times r_{i}}$, $\widehat{S}_{ij}\in\R{r_{i} \times r_{j}}$, $\widehat{V}_j\in\R{r_j \times B_j}$. Note that all low-rank blocks in a row $i$ have the same left factor $\widehat{U}_i$, and all the low-rank blocks in a column $j$ have the same right factor $\widehat{V}_j$. This is one of the defining features of the \HT matrix. We denote $r_i$ as the rank of $i$-th block row, excluding full-rank blocks, and $r_j$ is the rank of $j$-th block column, excluding full-rank blocks.

Let us define a block matrix $D\in\R{N\times N}$:
$$[D]_{ij} = \begin{cases}
\widehat{D}_{ij},          & \quad \text{if } A_{ij} \text{ is a close block} \\
0,                         & \quad \text{if } A_{ij} \text{ is a far block}
\end{cases}.$$ 
Note that $D$ is typically a block-sparse matrix; it contains a number of nonzero blocks per block row that is independent of dimension. 
Also, define a block matrix $F\in\R{N\times N}$:
$$[F]_{ij} = \begin{cases}
0,                     & \quad \text{if } A_{ij} \text{ is a close block} \\
\widehat{F}_{ij} = \widehat{U}_i \widehat{S}_{ij} \widehat{V}_j,                                             & \quad \text{if } A_{ij} \text{ is a far block}
\end{cases}.$$ 
The matrix $A$ is split into two matrices:
\begin{equation}
    \label{eq:h20}
    A = D + F,
\end{equation} 

To write equation~\eqref{eq:f} in matrix form, we define rectangular diagonal matrices $U_1\in \R{N\times N_2}$, $V_1\in \R{N_2\times N}$, where ${N_2=\sum_{i=1}^M r_i}$:
$$U_1 = \begin{bmatrix}
    \widehat{U}_{1} & & \\
    &\ddots & \\
    & & \widehat{U}_{M}
\end{bmatrix}, 
\quad V_1 = \begin{bmatrix}
    \widehat{V}_{1} & & \\
    &\ddots & \\
    & & \widehat{V}_{M}
\end{bmatrix},$$

We also define a matrix $S_1\in\R{N_2\times N_2}$ as:
$${\color{black}[S_1]_{ij}} = \begin{cases}
0,                & \quad \text{if } A_{ij} \text{ is a close block} \\
\widehat{S}_{ij}, & \quad \text{if } A_{ij} \text{ is a far block}
\end{cases}.$$ 






These definitions allow us to rewrite equation~\eqref{eq:f} in matrix form as:
$$F = U_1S_1V_1,$$
which allows us to express equation~\eqref{eq:h20} as:
$$A = D + U_1S_1V_1.$$





\begin{figure}[H]
\begin{center}
\begin{tikzpicture}[scale=0.34]
    \begin{scope}[shift={(0,0)}]
    \draw [draw=dgreen, fill=lgreen] (0, 0) rectangle (8, 8);
    \node at (4,-1) {$A$};
    \draw[step=1.0,dgreen,thin] (0,0) grid (8,8);
    \foreach \x in {0,...,7}{
        \draw [draw=lgreen,fill=dgreen] (\x, 7-\x) rectangle (\x+1, 8-\x);}
     \foreach \x in {0,...,6}{
        \draw [draw=lgreen,fill=dgreen] (\x, 6-\x) rectangle (\x+1, 7-\x);
        \draw [draw=lgreen,fill=dgreen] (\x+1, 7-\x) rectangle (\x+2, 8-\x);}
    \foreach \x in {0,...,2}{
        \draw [draw=lgreen,fill=dgreen] (3 + 2*\x, 7-2*\x) rectangle (4+2*\x, 8-2*\x);
        \draw [draw=lgreen,fill=dgreen] (2*\x, 4-2*\x) rectangle (1+2*\x, 5-2*\x);}

    \node at (10.,4) {\huge $\approx$};
    \end{scope}
    
    \begin{scope}[shift={(12,0)}]
    \draw [draw=dgreen] (0, 0) rectangle (8, 8);
        \node at (4,-1) {$D$};
        \draw[step=1.0,dgreen,thin] (0,0) grid (8,8);
        \node at (10.,4) {\huge $+$};
        \foreach \x in {0,...,7}{
        \draw [draw=lgreen,fill=dgreen] (\x, 7-\x) rectangle (\x+1, 8-\x);}
     \foreach \x in {0,...,6}{
        \draw [draw=lgreen,fill=dgreen] (\x, 6-\x) rectangle (\x+1, 7-\x);
        \draw [draw=lgreen,fill=dgreen] (\x+1, 7-\x) rectangle (\x+2, 8-\x);}
    \foreach \x in {0,...,2}{
        \draw [draw=lgreen,fill=dgreen] (3 + 2*\x, 7-2*\x) rectangle (4+2*\x, 8-2*\x);
        \draw [draw=lgreen,fill=dgreen] (2*\x, 4-2*\x) rectangle (1+2*\x, 5-2*\x);}
    \end{scope}
    
    \begin{scope}[shift={(24,0)}]
        \draw [draw=dgreen] (0, 0) rectangle (4, 8);
        \foreach \x in {0,...,7}{
            \draw [draw=dgreen,fill=dgreen] (\x*0.5, 7-\x) rectangle ((\x*0.5+0.5, 8-\x);}
        \node at (1,-1) {$U_1$};
    \end{scope}
    
    \begin{scope}[shift={(29,0)}]
        \draw [draw = lgreen,fill=lgreen] (0, 4) rectangle (4, 8);
        \draw[step=0.5,dgreen,thin] (0,4) grid (4,8);
        \foreach \x in {0,...,7}{
                \draw [draw=dgreen,fill=white] (\x*0.5, 7.5-\x*0.5) rectangle (\x*0.5+0.5, 8-\x*0.5);}
     \foreach \x in {0,...,6}{
        \draw [draw=dgreen,fill=white] (\x*0.5, 7-\x*0.5) rectangle (\x*0.5+0.5, 7.5-\x*0.5);
        \draw [draw=dgreen,fill=white] (\x*0.5+0.5, 7.5-\x*0.5) rectangle (\x*0.5+1, 8-\x*0.5);}
    \foreach \x in {0,...,2}{
        \draw [draw=dgreen,fill=white] (1.5 + \x, 7.5-\x) rectangle (2+\x, 8-\x);
        \draw [draw=dgreen,fill=white] (\x, 6-\x) rectangle (0.5+\x, 6.5-\x);}
    \node at (2,-1) {$S_1$};
    \end{scope}
    
    \begin{scope}[shift={(34,0)}]
        \draw [draw=dgreen] (0, 4) rectangle (8, 8);
        \foreach \x in {0,...,7}{
            \draw [draw=dgreen,fill=dgreen] (\x, 7.5-\x*0.5) rectangle (\x+1, 8-\x*0.5);
        }
        \node at (4,-1) {$V_1$};
    \end{scope}
\end{tikzpicture}
\end{center}
\vspace*{-12pt} 
\caption{Illustration of the one-level \HT matrix}
\label{fig:h21}
\end{figure}
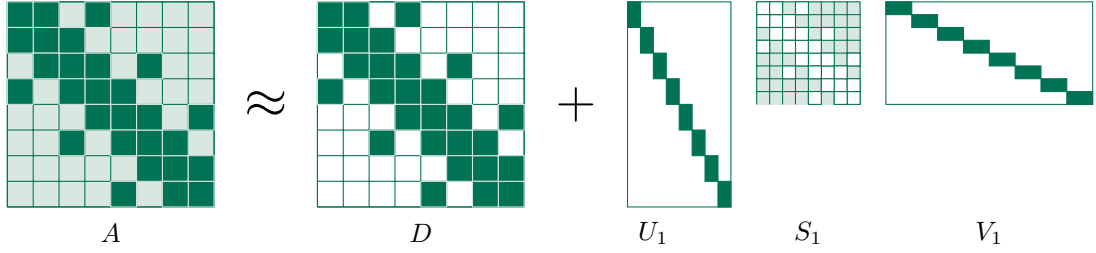
\begin{figure}[h]
\begin{center}
\tikzdeclarepattern{
  name=hatch,
  parameters={
      \pgfkeysvalueof{/pgf/pattern keys/size},
      \pgfkeysvalueof{/pgf/pattern keys/angle},
      \pgfkeysvalueof{/pgf/pattern keys/line width},
  },
  bounding box={
    (0,-0.5*\pgfkeysvalueof{/pgf/pattern keys/line width}) and
    (\pgfkeysvalueof{/pgf/pattern keys/size},
0.5*\pgfkeysvalueof{/pgf/pattern keys/line width})},
  tile size={(\pgfkeysvalueof{/pgf/pattern keys/size},
\pgfkeysvalueof{/pgf/pattern keys/size})},
  tile transformation={rotate=\pgfkeysvalueof{/pgf/pattern keys/angle}},
  defaults={
    size/.initial=5pt,
    angle/.initial=45,
    line width/.initial=.4pt,
  },
  code={
      \draw [line width=\pgfkeysvalueof{/pgf/pattern keys/line width}]
        (0,0) -- (\pgfkeysvalueof{/pgf/pattern keys/size},0);
  },
}

\begin{tikzpicture}[scale=0.2]
    \begin{scope}[shift={(0,0)}]
    \draw [draw=dgreen, fill=lgreen] (0, 0) rectangle (8, 8);
    \draw[step=1.0,dgreen,thin] (0,0) grid (8,8);
    \foreach \x in {0,...,7}{
        \draw [draw=dgreen,fill=white] (\x, 7-\x) rectangle (\x+1, 8-\x);}
     \foreach \x in {0,...,6}{
        \draw [draw=dgreen,fill=white] (\x, 6-\x) rectangle (\x+1, 7-\x);
        \draw [draw=dgreen,fill=white] (\x+1, 7-\x) rectangle (\x+2, 8-\x);}
    \foreach \x in {0,...,2}{
        \draw [draw=dgreen,fill=white] (3 + 2*\x, 7-2*\x) rectangle (4+2*\x, 8-2*\x);
        \draw [draw=dgreen,fill=white] (2*\x, 4-2*\x) rectangle (1+2*\x, 5-2*\x);}
    \node at (4,-2) {$S_1$};
    \node at (10.,4) {\large $=$};
    \end{scope}
    
    \begin{scope}[shift={(12,0)}]
    \draw [draw=dgreen, fill=lgreen] (0, 0) rectangle (8, 8);
    \draw[step=2.0,dgreen,thin] (0,0) grid (8,8);
    \foreach \x in {0,...,7}{
        \draw [draw=dgreen,fill=white] (\x, 7-\x) rectangle (\x+1, 8-\x);}
     \foreach \x in {0,...,6}{
        \draw [draw=dgreen,fill=white] (\x, 6-\x) rectangle (\x+1, 7-\x);
        \draw [draw=dgreen,fill=white] (\x+1, 7-\x) rectangle (\x+2, 8-\x);}
    \foreach \x in {0,...,2}{
        \draw [draw=dgreen,fill=white] (3 + 2*\x, 7-2*\x) rectangle (4+2*\x, 8-2*\x);
        \draw [draw=dgreen,fill=white] (2*\x, 4-2*\x) rectangle (1+2*\x, 5-2*\x);}
    \node at (4,-2) {$S_1$};
    \foreach \x in {0,...,3}{
    \draw[draw=dgreen,pattern={hatch[size=4pt,line width=1.8pt,angle=45]}, pattern color=dgreen, ] (\x*2,8-\x*2) rectangle (2+\x*2,6-\x*2);}
    \foreach \x in {0,...,2}{
    \draw[draw=dgreen,pattern={hatch[size=4pt,line width=1.8pt,angle=45]}, pattern color=dgreen, ] (\x*2,4-\x*2) rectangle (2+\x*2,6-\x*2);
    \draw[draw=dgreen,pattern={hatch[size=4pt,line width=1.8pt,angle=45]}, pattern color=dgreen, ] (\x*2+2,6-\x*2) rectangle (4+\x*2,8-\x*2);}
    \node at (10.,4) {\large $\approx$};
    \end{scope}
    
    \begin{scope}[shift={(24,0)}]
    \draw[step=2.0,dgreen,thin] (0,0) grid (8,8);
    \node at (4,-2) {$D_2$};
    \foreach \x in {0,...,3}{
    \draw[draw=lgreen, fill=dgreen] (\x*2,8-\x*2) rectangle (2+\x*2,6-\x*2);}
    \foreach \x in {0,...,2}{
    \draw[draw=lgreen, fill=dgreen] (\x*2,4-\x*2) rectangle (2+\x*2,6-\x*2);
    \draw[draw=lgreen, fill=dgreen] (\x*2+2,6-\x*2) rectangle (4+\x*2,8-\x*2);}
    \node at (10.,4) {\large $+$};
    \end{scope}
    
    \begin{scope}[shift={(36,0)}]
        \draw [draw=dgreen] (0, 0) rectangle (4, 8);
        \foreach \x in {0,...,3}{
            \draw [draw=dgreen,fill=dgreen] (\x, 6-\x*2) rectangle (\x+1, 8-\x*2);
            }
        \node at (2,-2) {$U_2$};
    \end{scope}
    
    \begin{scope}[shift={(42,0)}]
        \draw [draw = dgreen,fill=lgreen] (0, 4) rectangle (4, 8);
        \foreach \x in {0,...,3}{
                \draw [draw=dgreen,fill=white] (\x, 7-\x) rectangle (\x+1, 8-\x);}
        \foreach \x in {0,...,2}{
                \draw[draw=dgreen, fill=white] (\x,6-\x) rectangle (1+\x,7-\x);
                \draw[draw=dgreen, fill=white] (\x+1,7-\x) rectangle (2+\x,8-\x);}
        \draw[step=1.0, dgreen, thin] (0,4) grid (4,8);
        \node at (2,-2) {$S_2$};
    \end{scope}
    
    \begin{scope}[shift={(48,0)}]
        \draw [draw=dgreen] (0, 4) rectangle (8, 8);
        \foreach \x in {0,...,3}{
            \draw [draw=dgreen,fill=dgreen] (\x*2, 7-\x) rectangle (\x*2+2, 8-\x);
        }
        \node at (2,-1) {$V_2$};
    \end{scope}
\end{tikzpicture}
\end{center}
\vspace*{-12pt} 
\caption{Block low-rank factorization of the matrix $S_1$}
\label{fig:s1}
\end{figure}
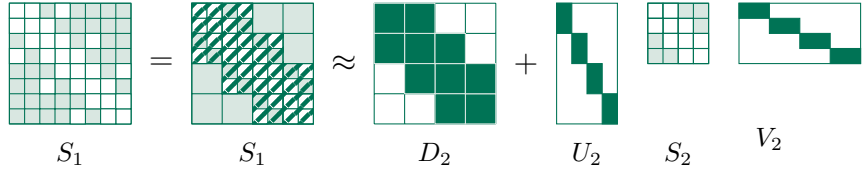

This decomposition is illustrated in Figure~\ref{fig:h21}. The figure illustrates a special case; the matrix $D$ may have more complex block structure. In general, it could be any block-sparse matrix. 
The matrix $S_1$ is a dense matrix. Our next goal is to compress it using the same strategy as we applied to matrix $A$.
We assemble together blocks of matrix $S_1$ in patches of size $p$. (Typical choices of $p$ are $2$ or $2^d$, where $d$ is the physical dimension of the problem.) We obtain matrix $S_1$ with $\frac{M}{p}$ block {\color{newblue} rows and $\frac{M}{p}$ block columns}. Those blocks again separate into a set of full-rank close and low-rank far blocks. Analogously to the previous procedure, we obtain the factorization
\begin{equation}
   S_1 =  D_2 + U_2 S_2 V_2, 
   \label{eq:s1}
\end{equation}
\noindent
with $D_2\in \R{N_2\times N_2}$, $U_2\in \R{N_2\times N_3}$, $V_2\in \R{N_3\times N_2}$, $S_2\in \R{N_3\times N_3}$. $N_3$ is $\sum_{i=1}^{\frac{M}{p}} r_i$, where $r_i$ are the ranks of the far blocks of the matrix $S_1$. An illustration of the factorization~\eqref{eq:s1} is shown in Figure~\ref{fig:s1}.


The process continues until the matrix $S_l$ on the $l^{th}$ level has low-rank blocks.
Assuming $l=2$, we obtain:
$$A = D +  U_1 \left(D_2 + U_2 S_2 V_2\right) V_1.$$
In the general case:
\begin{equation}
    A = D +  U_1 \left(D_2 + U_2 (\dots (D_l + U_lS_lV_l)\dots) V_2\right) V_1.
    \label{eq:h2}
\end{equation}
This recursive summation is called the \HT approximation of the matrix $A$.

\begin{figure}[H]
\begin{center}
\begin{tikzpicture}[scale=0.7]
    \begin{scope}[shift={(0,0)}]
        \def \lline {-7.5};
        \def \mline {-4};
        \def \uline {-0.5};
        \draw[-{Stealth[length=4mm]}, draw=dgreen, line width=0.5mm] (0,\lline)--(5.9,\lline);
        \draw[-{Stealth[length=4mm]}, draw=dgreen, line width=0.5mm] (1,\mline)--(4.9,\mline);
        \draw[-{Stealth[length=4mm]}, draw=dgreen, line width=0.5mm] (2,\uline)--(3.9,\uline);
        
        \draw[-{Stealth[length=4mm]}, draw=dgreen, line width=0.5mm](0,\lline)--(1,\mline);
        \draw[-{Stealth[length=4mm]}, draw=dgreen, line width=0.5mm](1,\mline)--(2,\uline);
        \path[fill=dgreen,draw=dgreen] (0,\lline) circle (1mm);
        \path[fill=dgreen,draw=dgreen] (1,\mline) circle (1mm);
        \path[fill=dgreen,draw=dgreen] (2,\uline) circle (1mm);
    
        \draw[-{Stealth[length=4mm]},draw=dgreen, line width=0.5mm](4,\uline)--(5,\mline);
        \draw[-{Stealth[length=4mm]},draw=dgreen, line width=0.5mm](5,\mline)--(6,\lline);
        \path[fill=dgreen,draw=dgreen] (4,\uline) circle (1mm);
        \path[fill=dgreen,draw=dgreen] (5,\mline) circle (1mm);
        \path[fill=dgreen,draw=dgreen] (6,\lline) circle (1mm);

        \draw (3,\lline-0.2) node [below]  {\large{$D$}} ;
        \draw (3,\mline-.2) node [below]  {\large{$D_2$}} ;
        \draw (3,\uline-0.2) node [below]  {\large{$S_2$}} ;
    
        \draw (0.5,\lline/2+\mline/2) node [left]  {\large{$V_1$}} ;
        \draw (1.5,\uline/2+\mline/2) node [left]  {\large{$V_2$}} ;
    
        \draw (5.5,\lline/2+\mline/2) node [right]  {\large{$U_1$}} ;
        \draw (4.5,\uline/2+\mline/2) node [right]  {\large{$U_2$}};
    
        \draw (-.1,\lline) node [left]  {\large{$x$}} ;
        \draw (.8,\mline+0.2) node [left]  {\large{$V_1x = x_1 $}} ;
        \draw (1.8,\uline+0.2) node [left] {\large{$V_2x_1 = x_2$}} ;
    
        \draw (6.1,\lline) node [right]  {\large{$y = Dx + U_1y_1$}} ;
        \draw (5.1,\mline+0.2) node [right]  {\large{$y_1 = S_1x_1+U_2y_2$}} ;
        \draw (4.1,\uline+0.2) node [right]  {\large{$y_2 = S_2 x_2$}} ;
    \end{scope}
    \end{tikzpicture}
\end{center}
\caption{Schematic of \HT matrix-by-vector multiplication, $l=2$}
\label{fig:h2_matvec}
\end{figure}
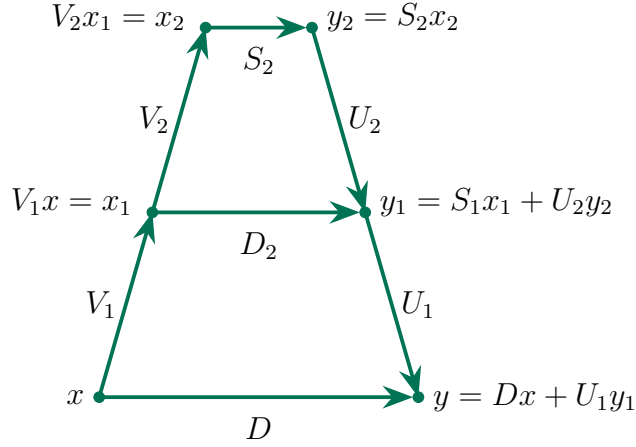

Multiplication of a vector $x\in \R{N}$ by the \HT matrix $A\in\R{N\times N}$  follows the formula~\eqref{eq:h2}:
$$y = Ax =  Dx +  U_1 \left(D_2 + U_2 (\dots (D_l + U_lS_lV_l)\dots) V_2\right)V_1x.$$
Letting $x_i = V_ix_{i-1}$, $x_i\in \R{N_i}$, $i = 1\dots l$, with $x = x_0$, we obtain:
$$y =  Dx +  U_1 \left(S_1x_1 + U_2 (\dots (S_{l-1}x_l + U_lS_lx_l)\dots)\right).$$
Letting $y_{i-1} = D_ix_{i-1} + U_iy_i$, $i = 1,\dots,l$,  $y_l = S_l x_l$, with $D_1 = D$ and $y_0 = y$, we obtain:
$$y =  Dx +  U_1y_1.$$
Figure~\ref{fig:h2_matvec} shows a schematic of the procedure.

\begin{remark}
\label{rem:h2on}
For the \HT matrix $A\in\R{N\times N}$, according to \cite{borm-h2-2010}, the storage requirement scales as $\OB(N)$. \HT construction and matrix-vector multiplication complexity is also $\OB(N)$. 
\end{remark}

\subsection {Multigrid method}

A multigrid method~\cite{bran-mg-1977,hack-mg-2013} is a powerful numerical technique used for solving partial differential equations (PDEs) and linear systems of equations with sparse matrices. While a solver in its own right, it is often used to precondition other iterative algorithms to help them converge rapidly.

The essential idea of multigrid is to define and solve problems on multiple grids of varying levels of resolution, referred to as the multigrid hierarchy. These grids range from coarse to fine, with each level representing a discretization of the problem at a different scale. At the coarsest level, the problem is smaller and easier to solve. Solutions obtained at this level are then recursively interpolated to finer grids. This hierarchical approach allows the method to address errors more effectively, each on an optimal scale, and accelerate convergence.

The multigrid method operates in cycles, typically consisting of phases of smoothing and inter grid transfer (restriction and prolongation). In ``smoothing'', iterative methods like Jacobi~\cite{jaco-jaco-1846}, Gauss-Seidel~\cite{seid-gs-1873}, CG~\cite{hest-cg-1952,shew-cg-1994}, etc., are applied to reduce errors. (Smoothing is in quotation marks, because this operation can sometimes ``roughen'' the solution.) In restriction, information is passed from finer to coarser grids and, in prolongation, from coarser to finer grids.

In this paper, as is often done in other contexts, we broaden the idea of multigrid beyond problems arising from discretized PDEs, here to systems with \HT matrices.
Analogously to using operators on multiple grid levels in a standard approach, we use the levels of \HT hierarchy, shown in Figure~\ref{fig:h2_matvec}. The transfer matrices $U_i$ and $V_i$, $i\in 1\dots l$ in \HT serve naturally as restriction and prolongation operators. The derivation of a multilevel iteration applied to the \HT system ($\mathcal{H}^{2}$-MG) is detailed in Section~\ref{sec:htmg}.

One of the key advantages of multigrid is its ability to rapidly converge to a highly accurate solution, often achieving convergence rates that are essentially independent of the problem size. This makes its translation to systems with \HT matrices highly beneficial and leads to a new and effective solver for systems with \HT matrices.

\subsection{$\mathcal{H}^{2}$-MG}
\label{sec:htmg}

In this section, we present a multigrid algorithm for a system with an \HT matrix.
Consider
\begin{equation}
    Ax = b,
    \label{eq:axb}
\end{equation}
\noindent
where $A\in\R{N\times N}$ is an $l$-level \HT matrix, $b\in \R{N}$ is the right-hand side, $x_0\in \R{N}$ is the initial guess for the solution. 

Let us follow the steps of the classic multigrid method, applying them to the system~\eqref{eq:axb}.
First, we reformulate the system in terms of error and residual. Let the residual of the system be
$$r_1 = b - Ax_0.$$
Then, subtracting from both parts of equation~\eqref{eq:axb} the $Ax_0$ term we write:
$$Ax - Ax_0 = b - Ax_0.$$
Noting that $Ax - Ax_0 = A(x-x_0) = Ae_1$, where $e_1 = x - x_0$ is the error, we obtain the error residual equation:
\begin{equation}
Ae_1 = r_1,    
\label{eq:aer}
\end{equation}
as an equivalent form of ~\eqref{eq:axb}.

The next stage of multigrid is the application of smoothing iterations to the system~\eqref{eq:axb}. 
We consider several iterations of some iterative method on the system~\eqref{eq:aer} as a smoother, $\textbf{Iter}(A,r_1,0)$. In  this notation, the first parameter is the operator, the second is the right-hand side, the third is the initial guess, and the result is the solution after the iterations.

In our numerical experiments, which are on symmetric positive definite matrices, we use a CG~\cite{shew-cg-1994, hest-cg-1952} solver as a smoother because of its fast and easily parallelized application. Thanks to the \HT structure of matrix $A$, the application of one iteration of CG is linear in time and memory cost.
After smoothing, we receive an approximation of the error:
$$\tilde{e}_1 = \textbf{Iter}(A,r_1,0).$$
We then compute the residual:
$$\hat{r}_1 = r_1 - A\tilde{e}_1$$
and subtract the term $A\tilde{e}_1$ from both parts of equation~\eqref{eq:aer} to obtain 
$$Ae_1 - A\tilde{e}_1= r_1 - A\tilde{e}_1.$$
Letting $\hat{e}_1 = e_1 - \tilde{e}_1$, we obtain the system
\begin{equation}
    A\hat{e}_1 = \hat{r}_1.
    \label{eq:aerhat}
\end{equation}
We then build the restriction and prolongation operators. 
If we write the matrix $A$, with orthogonal bases $U$ and $V$, explicitly in its \HT format, we obtain an expanded version of Equation~\eqref{eq:aerhat}:
$$ D +  U_1 \left(D_2 + U_2 (\dots (D_l + U_lS_lV_l)\dots) V_2\right) V_1 \hat{e}_1 = \hat{r}_1,$$
from which we see that a straightforward way to restrict the system is to multiply it by $U^{\top}_1$ on the left and to insert matrix $I = V^{\top}_1V_1$ between $A$ and $x$. We obtain the restricted system:
$$ U^{\top}_1(D +  U_1 \left(D_2 + U_2 (\dots (D_l + U_lS_lV_l)\dots) V_2\right) V_1 ) V^{\top}_1V_1\hat{e}_1  = U^{\top}_1\hat{r}_1 ,$$
or, if we open the first parentheses:
$$ (U^{\top}_1DV^{\top}_1 +  U^{\top}_1U_1 \left(D_2 + U_2 (\dots (D_l + U_lS_lV_l)\dots) V_2\right) V_1 V^{\top}_1)V_1 \hat{e}_1  = U^{\top}_1\hat{r}_1 .$$
Using the orthogonality of $U_1$ and $V_1$, this may be written as:
$$ (U^{\top}_1DV^{\top}_1 +  D_2 + U_2 (\dots (D_l + U_lS_lV_l)\dots) V_2)V_1 \hat{e}_1  = U^{\top}_1\hat{r}_1.$$
We define the restricted operator $A_2$ as:
\begin{equation}
    \label{eq:a2_str}
    A_2 = U^{\top}_1DV^{\top}_1 +  D_2 + U_2 (\dots (D_l + U_lS_lV_l)\dots) V_2.
\end{equation}
\begin{remark}
    Note that $A_2 \in \mathbb{R}^{N_2 \times N_2}$ has \HT structure, just like the matrix $A$. The matrices $A$ and $A_2$ have exactly the same structure, except that the matrix $A$ stores the full matrix $D$, while $A_2$ stores a reduced part of it, $U_1^{\top} D V_1^{\top}$. Thus, matrix $A_2$ stores less information than $A$. Since $A$ scales linearly as an \HT matrix, $A_2$ also scales linearly. See the detailed proof in Section~\ref{sec:compan}.
\end{remark}

We also defined the restricted error and residual vectors as:
$$e_2 = V_1 \hat{e}_1,$$
$$r_2 = U^{\top}_1\hat{r}_1.$$
In our case, the $U^{\top}_i$ matrices are the analogs of the restriction operators, $V^{\top}_i$ are the analogs of the prolongation operators, and the basis vectors of the \HT levels are the analogs of the coarser grids.  
We obtain the restricted system:
$$A_2e_2 = r_2.$$
Then, according to the multigrid algorithm, we apply a smoother to the restricted system and obtain an approximation of the error:
$$\tilde{e}_2 = \textbf{Iter}(A_2,r_2,0).$$
This continues until we reach the top level $l$. At this level, we have a system
$$A_l e_l = r_l,$$
\noindent
where $A_l$ is a small dense matrix since $l$ is the top level of \HT hierarchy. 
We solve the system directly. In our computations, we use the Cholesky factorization: 
$$e_l = \textbf{dir\_sol}(A_l, r_l).$$
Then, we move back from coarser to finer grids. 
We apply the prolongation operator $V_l^{\top}$ to the error $e_l$ and correct the $\tilde{e}_{l-1}$ error:
$$\tilde{e}_{l-1} = \tilde{e}_{l-1} + V_l^{\top} e_l .$$ 
Then, we apply the smoothing operator starting with the initial guess $\tilde{e}_{l-1}$:
$$e_{l-1} = \textbf{Iter}(A_{l-1},r_{l-1},{\color{black}\tilde{e}_{l-1}}).$$
We continue until we reach level 1. From the estimated error $e_1$, we obtain the approximation of the solution $x^*$:
$$x^* = x_0 + e_1.$$
This is analogous to the multigrid V-cycle. We can perform multiple V-cycles to obtain a more accurate solution, using $x^*$ as the initial guess for the next V-cycle. Figure~\ref{fig:h2mg1} is a visualization of \HTMG V-cycle. 

\begin{figure}
\begin{center}
\begin{tikzpicture}[scale=0.7]
    \def \lline {-7.5};
    \def \mline {-4};
    \def \uline {-0.5};
    \draw[-{Stealth[length=4mm]}, draw=dgreen, line width=0.5mm] (2,\uline)--(3.9,\uline);
    
    \draw[-{Stealth[length=4mm]}, draw=dgreen, line width=0.5mm](0,\lline)--(1,\mline);
    \draw[-{Stealth[length=4mm]}, draw=dgreen, line width=0.5mm](1,\mline)--(2,\uline);
    \path[fill=dgreen,draw=dgreen] (0,\lline) circle (1mm);
    \path[fill=dgreen,draw=dgreen] (1,\mline) circle (1mm);
    \path[fill=dgreen,draw=dgreen] (2,\uline) circle (1mm);

    \draw[-{Stealth[length=4mm]},draw=dgreen, line width=0.5mm](4,\uline)--(5,\mline);
    \draw[-{Stealth[length=4mm]},draw=dgreen, line width=0.5mm](5,\mline)--(6,\lline);
    \path[fill=dgreen,draw=dgreen] (4,\uline) circle (1mm);
    \path[fill=dgreen,draw=dgreen] (5,\mline) circle (1mm);
    \path[fill=dgreen,draw=dgreen] (6,\lline) circle (1mm);


    \draw (0,\lline-.1) node [below]  {\large{$\tilde{e}_1=\textbf{Iter}(A,{\color{newblue}r_1},0)$}} ;
    \draw (6.1,\lline-.1) node [below]  {\large{$e_1 = \textbf{Iter}(A, r_1, \tilde{e}_1)$}} ;

    \draw (0.5,\lline/2+\mline/2) node [left]  {\large{$r_2= U_1^{\top}(r_1-A\tilde{e}_1)$}} ;
    \draw (1.5,\uline/2+\mline/2) node [left]  {\large{$r_3= U_2^{\top}(r_2-A_2\tilde{e}_2)$}} ;

    \draw (5.5,\lline/2+\mline/2) node [right]  {\large{$\tilde{e}_1 = \tilde{e}_1 + V_1^{\top}e_2$}} ;
    \draw (4.5,\uline/2+\mline/2) node [right]  {\large{$\tilde{e}_2 = \tilde{e}_2 + V_2^{\top}e_3$}};

    \draw (.8,\mline+.2) node [left]  {\large{$\tilde{e}_2=\textbf{Iter}(A_2,r_2,0)$}} ;
    \draw (5.1,\mline+.2) node [right]  {\large{$e_2 = \textbf{Iter}(A_2, r_2, \tilde{e}_2) $}};
    
    \draw (3,\uline+.5) node [above]  {\large{$e_3 = \textbf{dir\_sol}(A_3,r_3)$}};
    \end{tikzpicture}
\end{center}
\caption{Schematic of the \HTMG algorithm}
\label{fig:h2mg1}
\end{figure}

\begin{remark}
    The visualization of \HTMG V-cycle, presented in Figure~\ref{fig:h2mg1} emphasizes the analogy of the \HTMG method with the \HT structure. Compare \HTMG V-cycle in Figure~\ref{fig:h2mg1} and the \HT matrix-vector product in Figure~\ref{fig:h2_matvec}. 
\end{remark}

In Algorithm~\ref{alg:htmg}, we give the formal description of one V-cycle of the \HTMG algorithm.  $A = A_1$, matrices $A_{i}\in\R{N_i\times N_i}$, $i=1,\dots, l-1$, are the \HT matrices, matrix $A_l\in\R{N_l\times N_l}$ is dense, $x_0\in \R{N}$ is the initial guess, and $b\in \R{N}$ is a right-hand side.

\begin{algorithm}[H]
\caption{One V-cycle of the \HTMG algorithm}
\label{alg:htmg}
\begin{algorithmic}[1]
\State $r_1 = b - A_1x_0$
\For{${\color{black}i} = \texttt{1\dots (}l\texttt{-1)}$}
    \State $\tilde{e}_i = \textbf{Iter}(A_i,r_i,0)$
    \State $r_{i+1} = U^{\top}_{i}(r_i - A_i\tilde{e}_i)$
\EndFor
\State $e_l = \textbf{dir\_sol}(A_l,r_l)$
\For{${\color{black}i} \texttt{ = (}l\texttt{-1)\dots 1}$}
    \State $\tilde{e}_i = \tilde{e}_i + V^{\top}_{i}e_{i+1}$
    \State $e_i = \textbf{Iter}(A_i,r_i,\tilde{e}_i)$
\EndFor
\State $x^* = x_0 + e_1$
\end{algorithmic}
\end{algorithm}
\noindent
We can run successive V-cycles using the output $x^*$ as the next initial guess. It may be of interest to consider multigrid cycles beyond V-cycles, namely W- or F-cycles.

\begin{remark}
    \label{rem:fine_coarse}
   In \HTMG, we treat the number of smoothing iterations on the finest level (with matrix $A$) and on the coarser levels (with matrices $A_i$, $i=2,\dots,l$) as two separate parameters. We denote the number of fine-level iterations as $n_f$ and the number of coarse-level iterations as $n_c$. Unlike the standard multigrid method, in our case, we have a physical grid only at the finest level, while all other grids are ``basis'' grids of the \HT matrix. Therefore, we treat the finest level of the problem differently by assigning it an independent number of smoothing iterations. In the numerical section, we empirically demonstrate that this approach leads to better convergence.
\end{remark}

\subsection{Complexity analysis}
\label{sec:compan}
In this subsection, we describe the time and memory complexity of the \HTMG algorithm.

The crucial feature of the \HTMG complexity analysis is the linear complexity of the \HT matrix. According to Remark~\ref{rem:h2on}, for the \HT matrix of size $N\times N$, the construction complexity, the memory requirements, and matrix-by-vector product complexity are all $\OB(N)$. Assume $c_{H2}$ to be the \HT matrix-vector product constant.

We first compute the \HTMG computational complexity  $n_{\text{op}}$ to run one V-cycle. The complexity of the fine grid smoothing iterations is $n_fc_{H2}N$, the coarse grid smoothing iterations is $n_cc_{H2}N_i$, for $i=2,\dots, l-1$, the complexity of the direct solver is $c_dN_l^3$, where $c_d$ is the direct solver complexity constant, and the restriction and prolongation operator complexity is negligible, compared to the 
smoothing iterations. The overall complexity is:
$$n_{\text{op}} = 2n_fc_{H2}N+ 2n_c\sum_{i=2}^{l-1} c_{H2}N_i+c_dN_l^3.$$
Assume, for simplicity, that the block size is $B$ for all blocks and the block rank $r$ is fixed for all levels (consider it to be the maximum rank of any block). Also, assume the number of blocks of the initial matrix $A$ is $M$, and the number of blocks stacked together while transferring from level to level is $d$. Thus, $N=MB$, $N_i = \frac{Mr}{d^{i-2}}$, for $i=2,\dots, l$, and the overall complexity is:
$$n_{\text{op}} = 2n_fc_{H2}MB+ 2n_c\sum_{i=2}^{l-1} \frac{Mrc_{H2}}{d^{i-2}}+c_dN_l^3.$$
Using the sum of a geometric progression, we obtain:
$$n_{\text{op}} = 2n_fc_{H2}MB+ 2n_c \left(\frac{d-\frac{1}{d^{l-3}}}{d-1}\right)c_{H2}Mr+c_{d}N_l^3.$$
By the construction of $\mathcal{H}^{2}$, the size of the coarsest level $N_l$ is a constant; thus,
$$n_{\text{op}} = \OB(N),$$
with the constant
$$c_{\text{op}} = 2n_fc_{H2}+ 2n_c \left(\frac{d-\frac{1}{d^{l-3}}}{d-1}\right)c_{H2}\frac{r}{B}$$
We next compute the memory requirements. The \HTMG algorithm requires storage of matrices $A$, $A_i$ for $i=2,\dots, l$, $U_i$, and $V_i$ for $i=1,\dots, l-1$. Matrices $U_i$ and $V_i$ are block-diagonal, and their storage is negligible compared to matrices $A$ and $A_i$. Assume \HT memory constant to be $c_m$. Thus, the memory requirements are:
$$n_{\text{mem}} = c_mN+ \sum_{i=2}^{l-1} c_mN_{i}+N_l^2.$$
Analogously to the time complexity, we rewrite matrix sizes in terms of $M$, $B$, $r$, and $d$ and sum the geometric progression to obtain:

$$n_{\text{mem}} = c_mN+ \left(\frac{d-\frac{1}{d^{l-3}}}{d-1}\right)c_mMr+N_l^2.$$ Taking into account that $N_l$ is a constant, we obtain:
$$n_{\text{mem}} = \OB(N),$$
with the overall \HTMG storage constant
$$c_{\text{mem}} = c_m+ \left(\frac{d-\frac{1}{d^{l-3}}}{d-1}\right)c_m\frac{r}{B}.$$
Thus, one V-cycle of the \HTMG algorithm is linear in both time and memory, though its constant factor is larger compared to a single iteration of the CG algorithm. However, in practice, we observe that the number of \HTMG V-cycles needed to achieve the required accuracy is independent of problem size, whereas the number of CG iterations needed to reach the same accuracy increases with the problem size.



\section{Numerical results}
The numerical experiments showcase our Python implementation of the \HTMG algorithm, demonstrating \HTMG convergence for two different kernels {\color{black}and a boundary element method (BEM) example}. The matrices are cast into the \HT format using the MCBH~\cite{mikh-mcbh-2016,mikh-mcbh-2018} method. The code for the \HTMG algorithm is publicly available at \href{https://github.com/dsushnikova/h2mg}{\url{https://github.com/dsushnikova/h2mg}}.

We compare the \HTMG algorithm with solvers CG and FMM-LU~\cite{sush-fmmlu-2023}. {\color{black}All solvers are implemented in Python without parallelism, and the experiments are run on a MacBook Pro (Apple M1 Pro, 16 GB RAM).} All three methods benefit from the \HT structure of matrix $A$; thus, the comparison is natural. The \HTMG algorithm uses CG iteration for smoothing, illustrating how the coarse-level iterations improve the convergence relative to CG alone. 

{\color{newblue} In the kernel matrix tests}, we randomly generate $x_{\text{true}}$, then compute $b = A x_{\text{true}}$ and solve the system $A x = b$ for~$x$. {\color{black} Then we plot convergence of the $A$-norm of the error $||e_k||_A$ vs.\ iteration count $k$, where the $A$-norm refers to the energy norm induced by matrix $A$.}

\begin{remark}
    The right-hand sides were chosen in this manner so that the error norm at each iteration could be easily computed.  We also tested a few cases where the right-hand sides were chosen randomly from a standard Gaussian distribution.  In these cases, the linear systems were harder to solve, but \HTMG outperformed unpreconditioned CG, just as in the detailed results that we show below.
\end{remark}

\subsection{Gaussian kernel}
\label{seq:num_gau}
For a first example, we consider the system
\begin{equation}
    \label{eq:axb_num}
    Ax = b,
\end{equation} with the Gaussian matrix $A$. We
consider a uniform tensor grid on a unit square $P\subset \R{2}$: $p_i\in P$, $i\in 1\dots N$,  where $N$ is the number of points.
The kernel matrix $A$ is given by the formula:
\begin{equation}
    \label{eq:a_gau}
    a_{ij} = \begin{cases}
    \exp{(-\frac{|p_i-p_j|^2}{\sigma})},          & \quad \text{if } i\neq j\\
    1 + c,                                        & \quad \text{if } i=j
    \end{cases},
\end{equation}

\noindent
where $\sigma\in\R{}$ is the dispersion parameter, $c\in\R{}$ is a small regularization parameter. 
Matrix $A$ is approximated {\color{newblue} in} the \HT format with accuracy $\epsilon=10^{-9}$ ($\epsilon$ is a relative error between the result of a matrix-vector multiplication performed using the \HT matrix and that using the original matrix), and the number of levels is chosen adaptively.

\subsubsection{Gaussian matrix, analysis of error evolution}
In our tests, we randomly generate $x_{\text{true}}$, then compute $b = A x_{\text{true}}$ and solve the system $A x = b$ for $x$. During the process, we can compute the error $e_{1*} = x_k - x_{\text{true}}$ since we know $x_{\text{true}}$. In the classic multigrid algorithm, the error $e_{1*}$ should behave in the following way: for the fine grid iterations, the high-frequency components of $e_{1*}$ should decrease rapidly, and during the coarse grid iterations, successive lower frequency components should decrease in turn.

In this section, we study the behavior of the error $e_{1*}$ in the \HTMG method.
Consider $U_1$, the transfer matrix from the finest level to the coarser one. $U_1$ has orthonormal columns. Assume that the matrix $Q_1$ has columns that span the complementary space.

To analyze the components of the error after smoothing in the basis of the interpolation matrix $U_1$, or in other words, the part of the error that will be projected to the coarser levels, we compute {\color{newblue}$(U_1)^{\top}e_{1*}$.} 
Similarly, for the part of the after-smoothing error in the basis of the matrix $Q_1$, {\color{newblue}$(Q_1)^{\top} e_{1*}$.}

\begin{figure}[H]
     \centering
     \subcaptionbox{{\color{newblue}Vector components of $(U_{1})^{\top} e_{1*}$\\
     (post-smoothing error in the basis of  $U_1$)}\label{fig:gau_err_uq1}}{
        \includegraphics[width=0.45\textwidth]{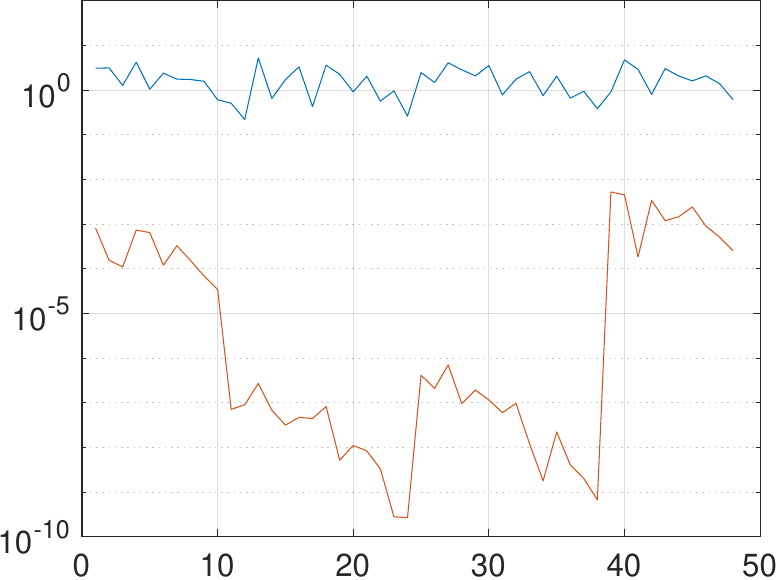}}
     \subcaptionbox{{\color{newblue}Vector components of $(Q_{1})^{\top} e_{1*}$\\
     (post-smoothing error in the basis of $Q_1$)} \label{fig:gau_err_uq2}}{
        \includegraphics[width=0.45\textwidth]{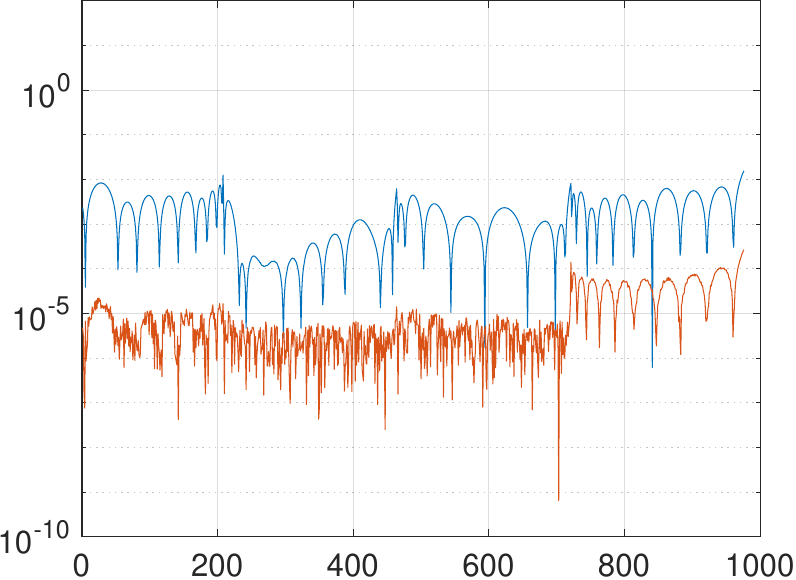}}
    \caption{Convergence comparison for {\color{black}two successive} coarse iterations for the \HTMG method}
    \label{fig:gau_err_uq}
\end{figure}

We study the 1D Gaussian matrix with points equally spaced between $0$ and $1$, $N=1024$, $\sigma = 0.01$, $c=10^{-3}$. The results of this analysis are presented in Figure~\ref{fig:gau_err_uq}. The blue curves represent the error components after the first smoothing step in the initial $V$-cycle, the red curves are for the second V-cycle. 

 The graphs clearly shows that the error component {\color{newblue}$(U_{1})^{\top} e_{1*}$} is larger than that in the error component {\color{newblue}$(Q_{1})^{\top} e_{1*}$} components, {\color{black}demonstrating} that the CG smoothing can complement coarse grid correction effectively.

\subsubsection{Gaussian kernel matrix, analysis of number of smoothing steps}
Consider the system \eqref{eq:axb_num} with matrix $A$ given by~\eqref{eq:a_gau} with $c = 10^{-3}$, $\sigma = 0.1$. In Remark~\ref{rem:fine_coarse}, we considered assigning the finest level of the problem a different number of smoothing steps than the coarser levels. In this subsection, we analyze the effect of parameters $n_f$ and $n_c$ on the \HTMG convergence.

\begin{figure}[t]
     \centering
     \subcaptionbox{{\color{black}$A$-norm of the error per coarse grid iteration} \label{fig:gau_fine_iter}}{
        \includegraphics[width=0.485\textwidth]{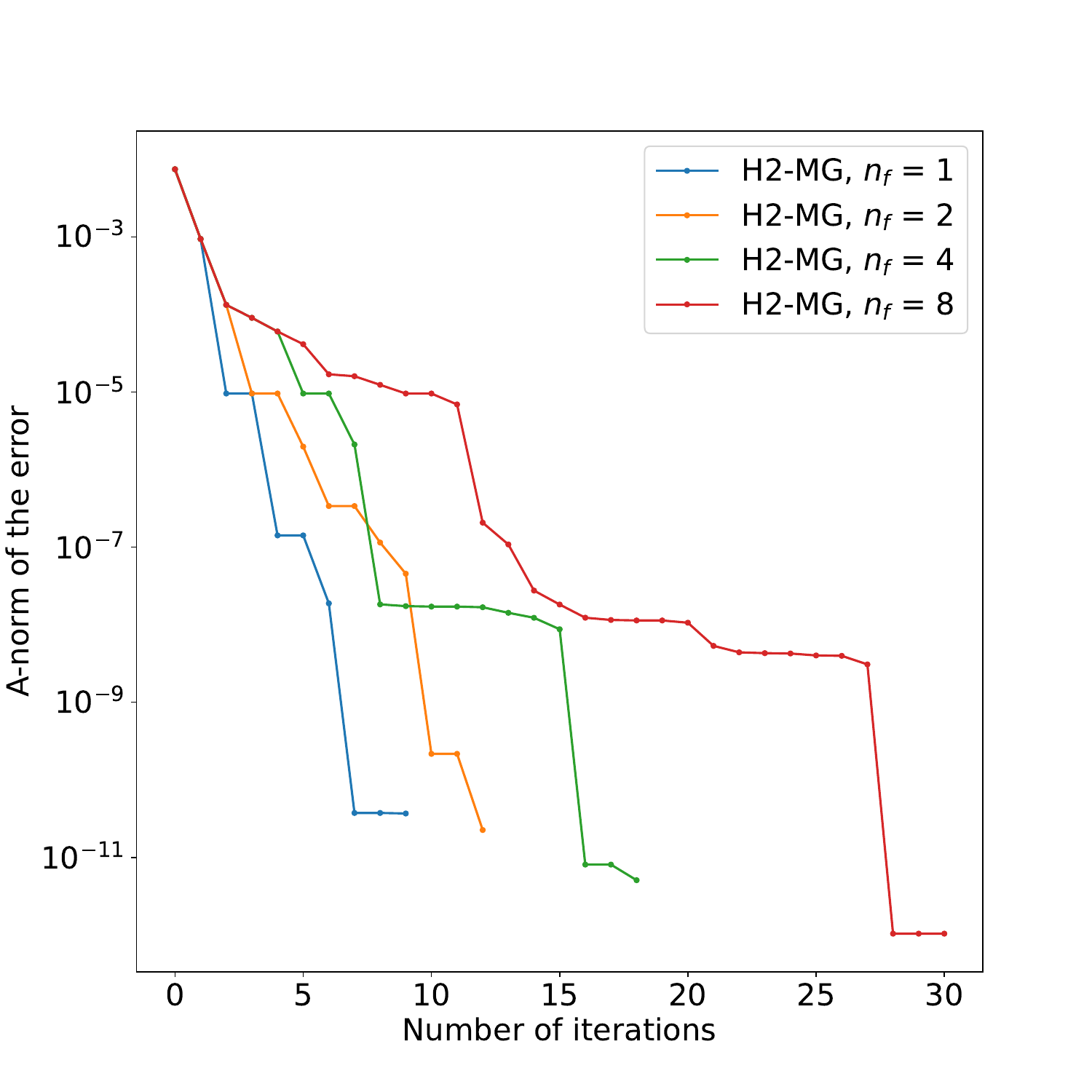}}
     \subcaptionbox{{\color{black}$A$-norm of the error per time}\label{fig:gau_fine_time}}{
        \includegraphics[width=0.485\textwidth]{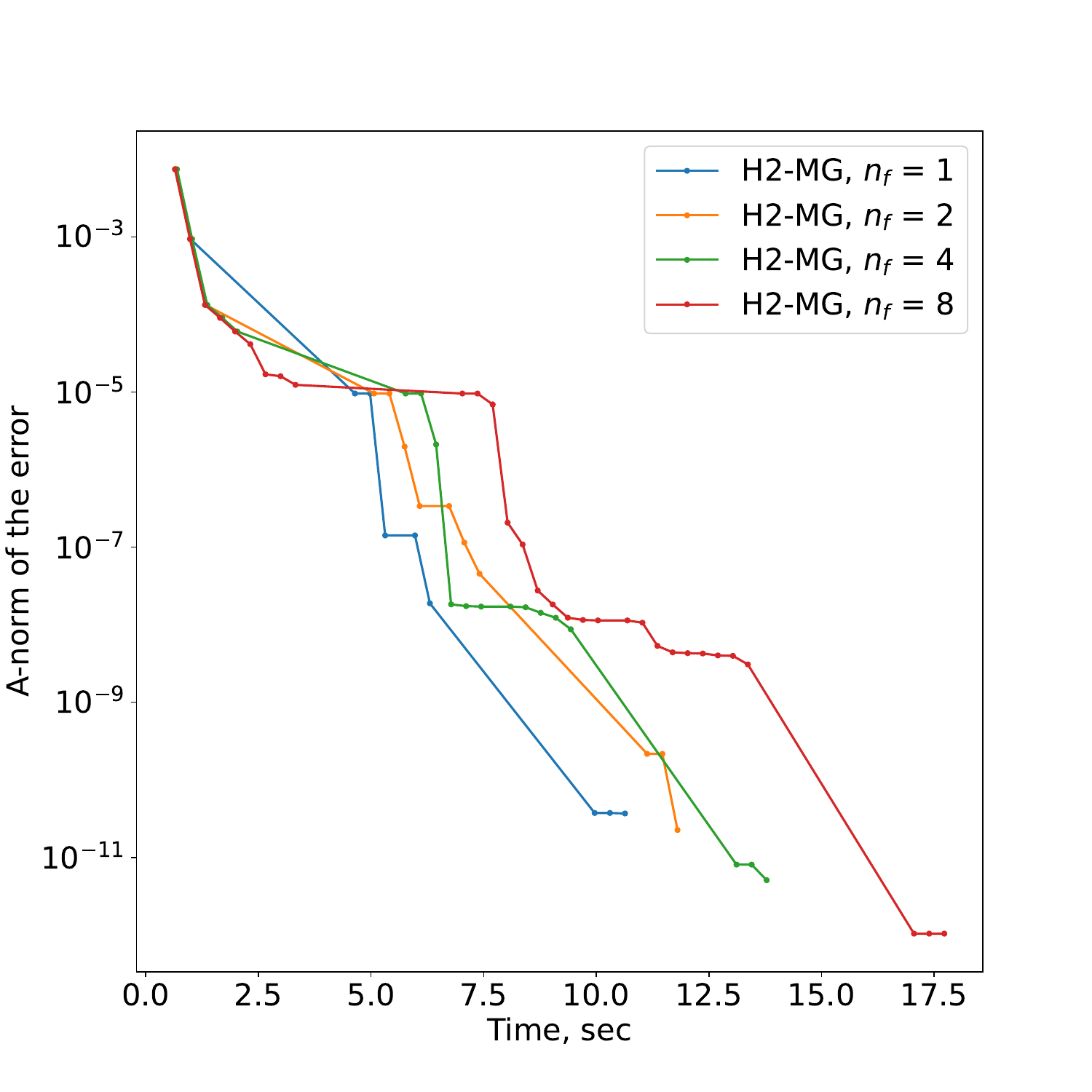}}
    \caption{Convergence comparison 
    for different numbers of fine iterations for the \HTMG method}
    \label{fig:gau_fine}
\end{figure}

We first study the convergence rate of the \HTMG depending on the number of fine iterations. Figure~\ref{fig:gau_fine} presents the convergence comparison for different numbers of fine iterations, $n_f$ to the accuracy $\varepsilon = 10^{-9}$, where $\varepsilon = \frac{(eAe^{\top})^{\frac{1}{2}}}{||b||_2} < 10^{-9}$. The number of coarse iterations is fixed $n_c = 40$. Fine and coarse iteration comparisons are given for matrix size $N=8\times 10^4$.
{\color{black} For \HTMG, we can track only the $A$-norm of the error on the finest level, which makes it especially useful for comparison with CG. In all figures below, we plot this error per finest level iterations (outer iterations). Since coarse-level iterations are not shown directly in these plots, we also include error-per-time plots to reflect the computational effort spent on coarser levels.}

{\color{newblue}We can see that using a smaller $n_f$ parameter leads to faster $A$-norm convergence, both when measured per outer iteration and per time.}

\begin{figure}
     \centering
     \subcaptionbox{{\color{black}$A$-norm of the error per coarse grid iteration}\label{fig:gau_coar_iter}}{
        \includegraphics[width=0.48\textwidth]{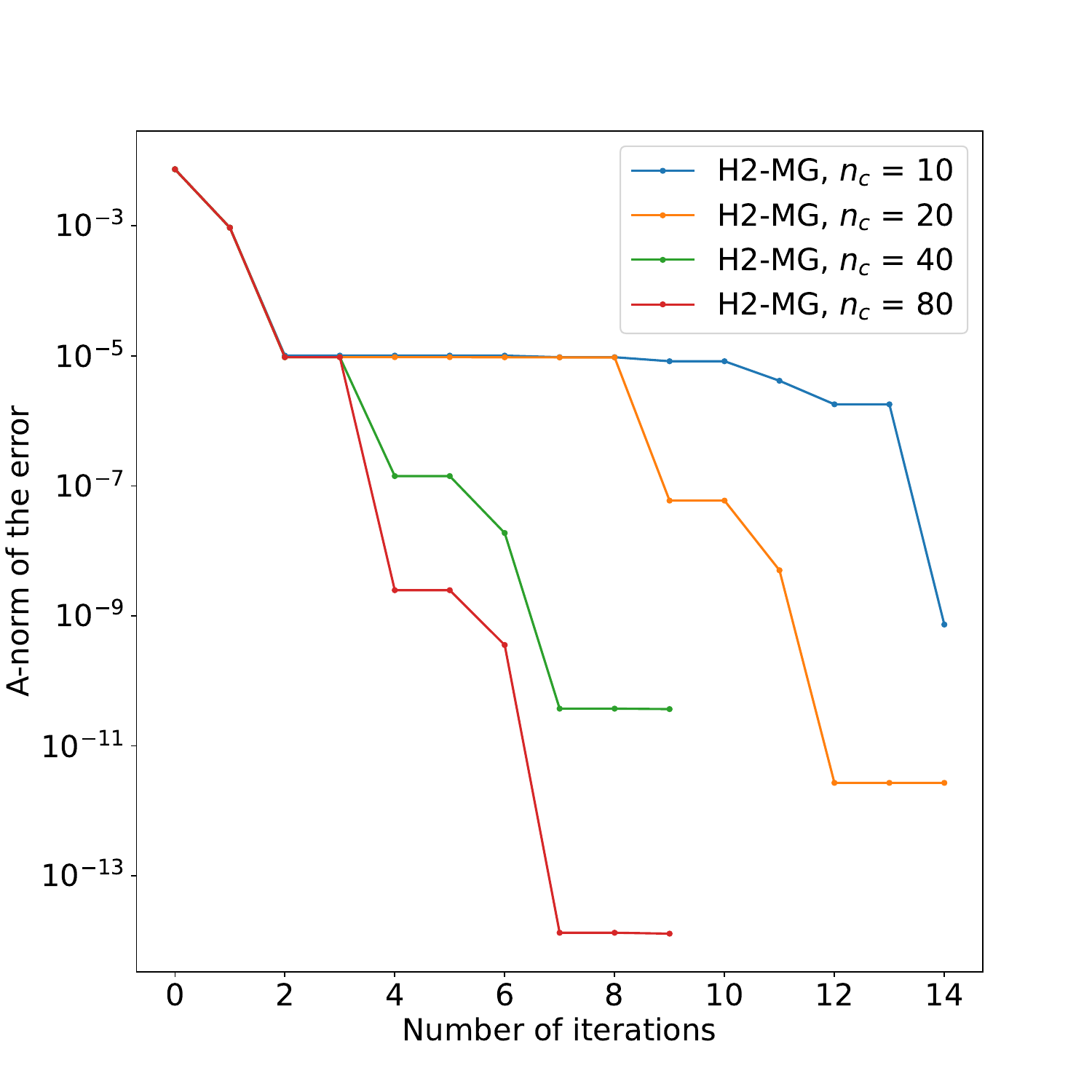}}
     \subcaptionbox{{\color{black}$A$-norm of the error per time} \label{fig:gau_coar_time}}{
        \includegraphics[width=0.48\textwidth]{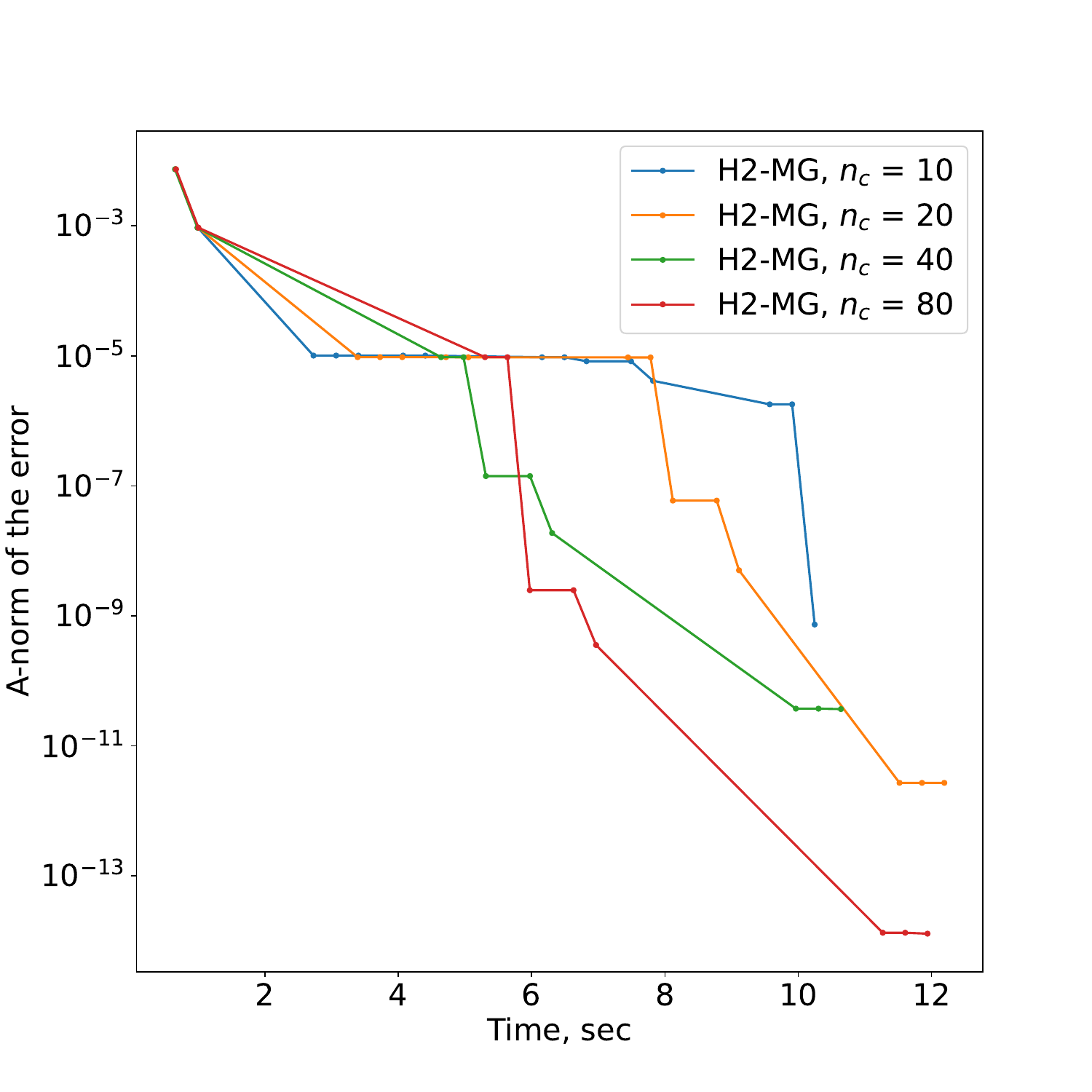}}
    \caption{Convergence comparison for different numbers of coarse iterations for the \HTMG method}
    \label{fig:gau_coarse}
\end{figure}

Figure~\ref{fig:gau_coarse} presents the convergence comparison to the accuracy $\varepsilon = 10^{-9}$ for different numbers of coarse iterations, $n_c$. The number of fine iterations is fixed at $n_f = 1$.

{\color{newblue} In this case, the $n_c$ parameter directly influences the speed of $A$-norm convergence per outer iteration: the more $n_c$, the faster the convergence per outer iteration. However, in terms of convergence per time,} the optimal number of coarse iterations is an intermediate value ($n_c = 40$). This is because there is a trade-off between the time spent on additional coarse iterations and the resulting improvement in convergence speed.

\begin{remark}
    Unlike a traditional multigrid method, where all levels are similar and require the same number of smoothing iterations, \HTMG uses a different number of smoothing iterations for the fine and coarse grids. The finest level is unique, as it is based on the physical grid, while all other levels are basis-induced, requiring a different number of smoothing steps.
\end{remark}

\subsubsection{Gaussian kernel matrix, asymptotics analysis}

\begin{figure}[H]
     \centering
     \subcaptionbox{ $c = 10^{-3}$ \label{fig:gau_nvc_s01_34}}{
        \includegraphics[width=0.48\textwidth]{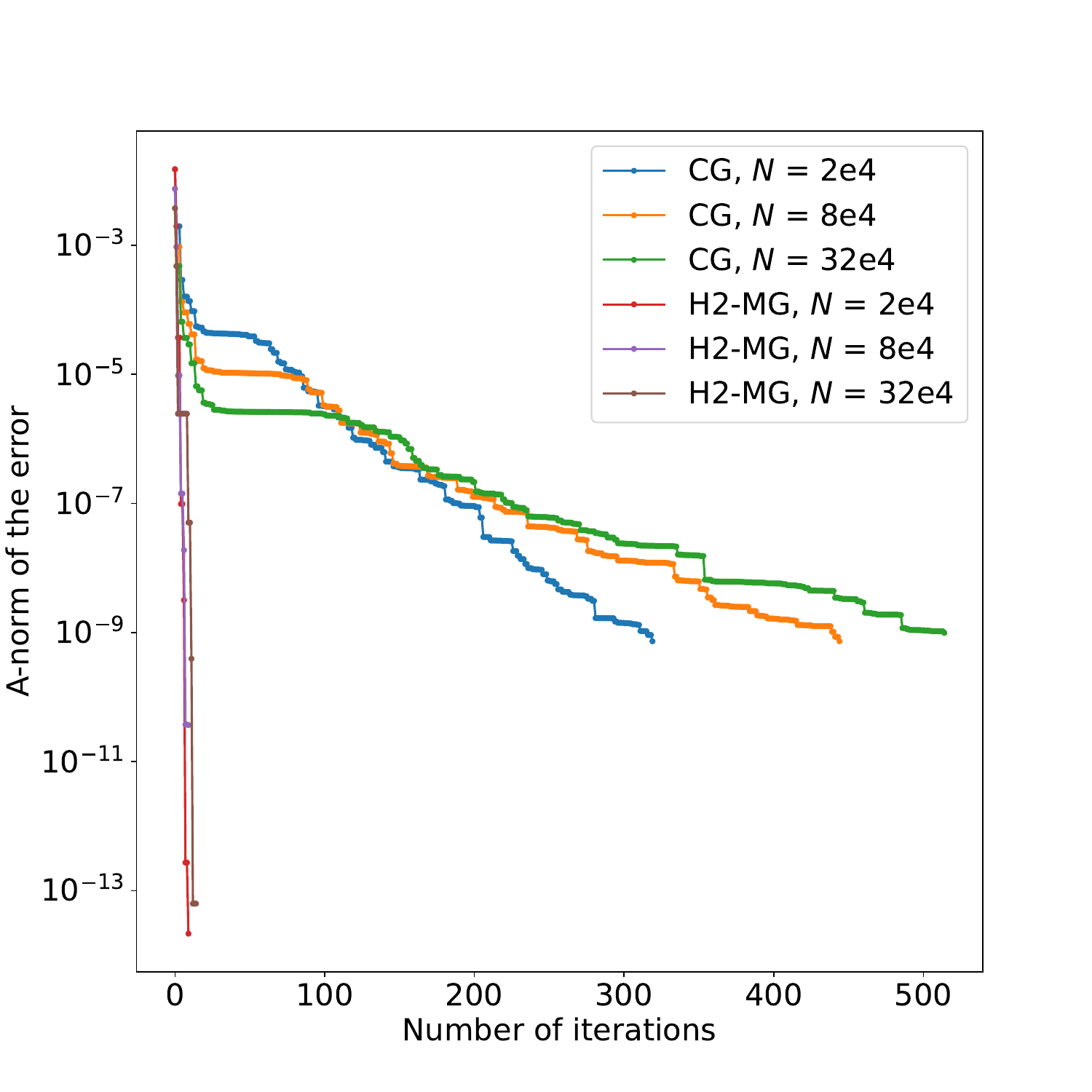}}
     \subcaptionbox{$c = 10^{-5}$ \label{fig:gau_nvc_s01_35}}{
        \includegraphics[width=0.48\textwidth]{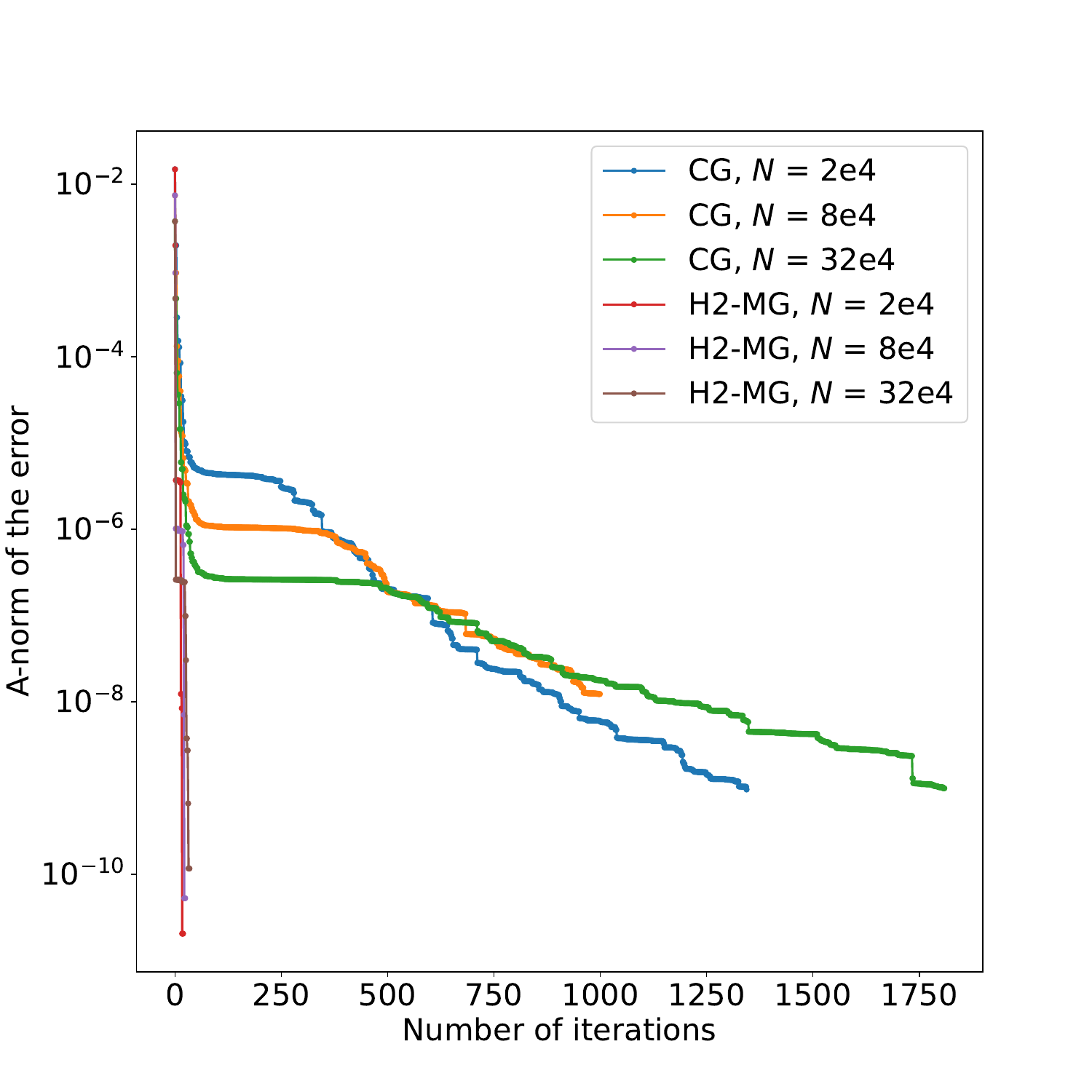}}
        \caption{Convergence evolution of \HTMG and CG as problem size increases,  $\sigma = 0.1$}
        \label{fig:gau_nvc_s01}
\end{figure}

\begin{figure}[H]
    \subcaptionbox{ $c = 10^{-3}$ \label{fig:gau_nvc_s01_36}}{
        \includegraphics[width=0.49\textwidth]{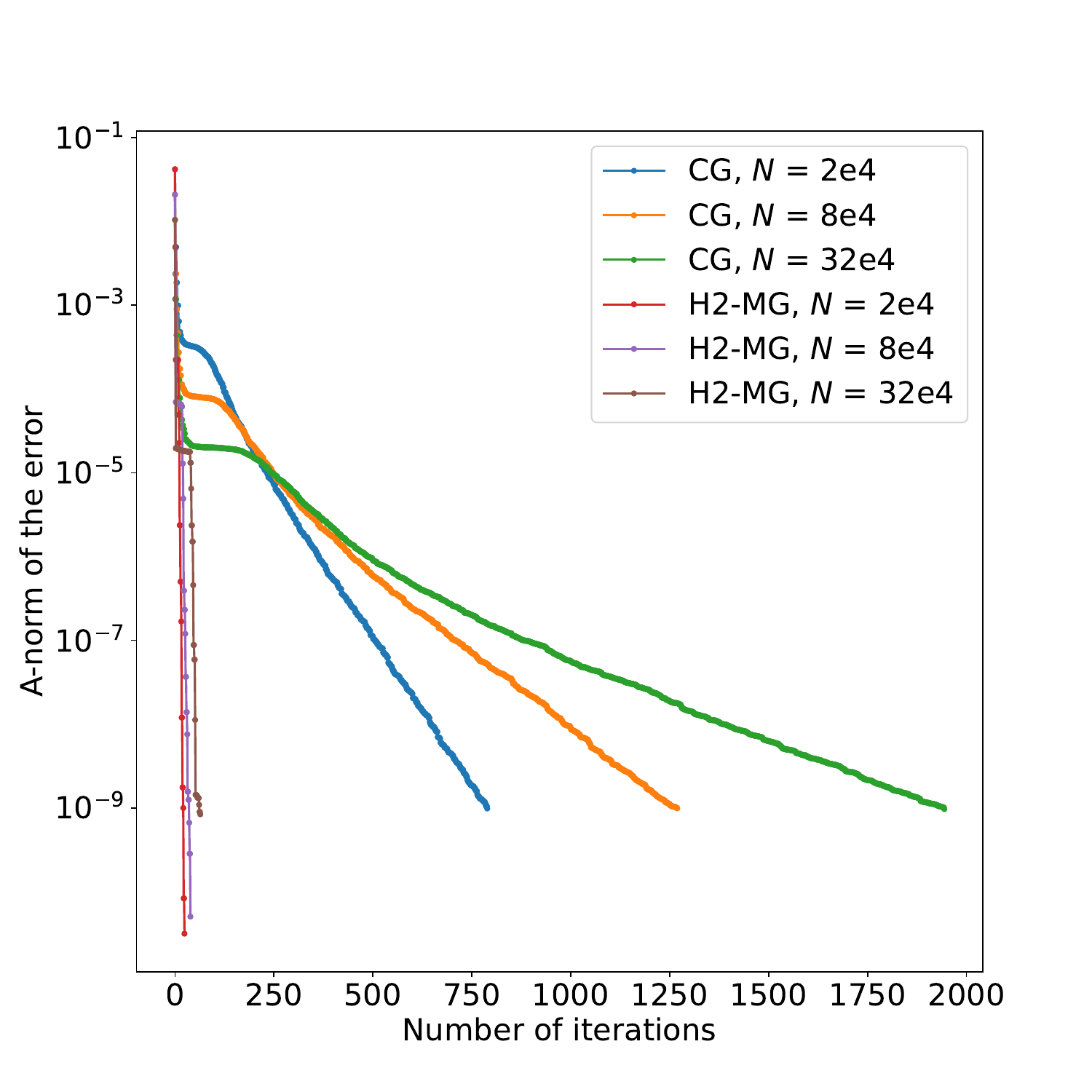}}
    \subcaptionbox{ $c = 10^{-5}$ \label{fig:gau_nvc_s01_37}}{
        \includegraphics[width=0.49\textwidth]{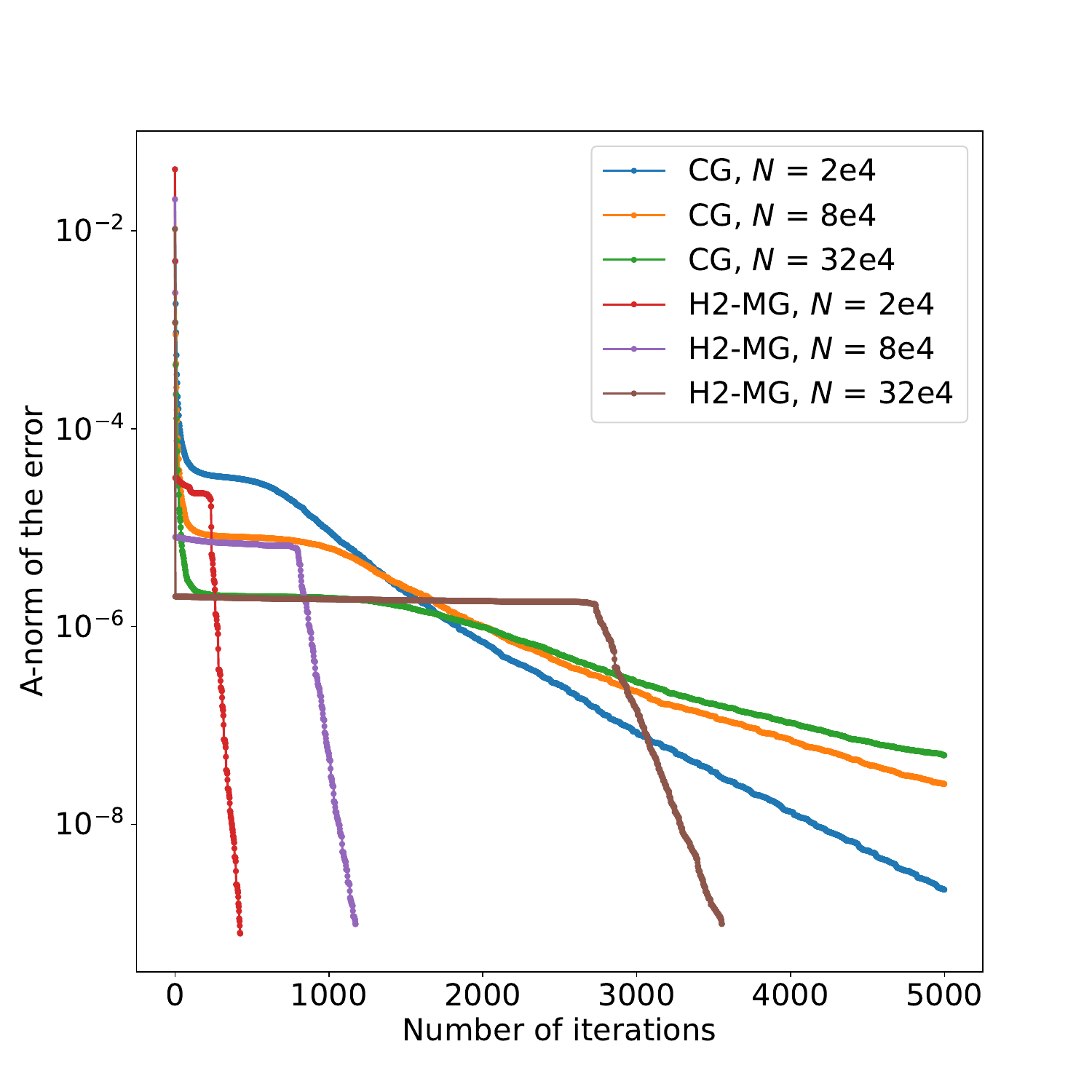}}
        \caption{Convergence evolution of \HTMG and CG as problem size increases,  $\sigma = 0.01$}
        \label{fig:gau_nvc_s001}
\end{figure}

In this subsection we consider the asymptotic behavior of the system~\eqref{eq:axb_num} with matrix $A$ given by~\eqref{eq:a_gau} for various combinations of parameters $c$ and $\sigma$.
Figures~\ref{fig:gau_nvc_s01}~and~\ref{fig:gau_nvc_s001} show the convergence of the \HTMG algorithm compared to CG across various problem sizes. 
We use a tolerance $\varepsilon = 10^{-9}$ and 
the parameters $n_f = 1$, $n_c = 40$.

We observe that the strong regularization parameter $c=10^{-3}$ leads to the fast convergence of the \HTMG method for both matrices with $\sigma = 0.1$ and $\sigma = 0.01$, while the weak regularization parameter $c=10^{-5}$ leads to the divergence of both algorithms for the larger problem sizes. 

The numbers of V-cycles required for convergence to a fixed accuracy are presented in Table~\ref{tab:gau_d2_vc}.

\begin{table}[h]
\centering
\begin{tabular}{c||c|c|c|c|c|c}
    \multirow{2}{*}{Matrix Parameters} & \multicolumn{6}{c}{Problem Size} \\ \cline{2-7}
    & 1e4 &  2e4   &  4e4 & 8e4 &  16e4 &  32e4  \\ \hline
   $\sigma = 0.1$, $c = 10^{-3}$  & 2 &  2    &  2  & 2   &  2   &  3     \\
   $\sigma = 0.1$, $c = 10^{-5}$  & 4 &  4    &  3  & 5  &  4   &  7     \\
   $\sigma = 0.01$, $c = 10^{-3}$  & 7 &  5   &  7  & 8   &  11  &  12    \\
   $\sigma = 0.01$, $c = 10^{-5}$  & 200 &  85   &  156 &  235  &  428  &  711    \\
\end{tabular}
\caption{Number of V-cycles to solve the system for different problem sizes}
\label{tab:gau_d2_vc}
\end{table}


As in the previous example, the number of V-cycles required for convergence to the fixed accuracy does not grow significantly.

\begin{figure}[H]
     \centering
     \subcaptionbox{ $c = 10^{-3}$ \label{fig:gau_st_s01_3}}{
        \includegraphics[width=0.48\textwidth]{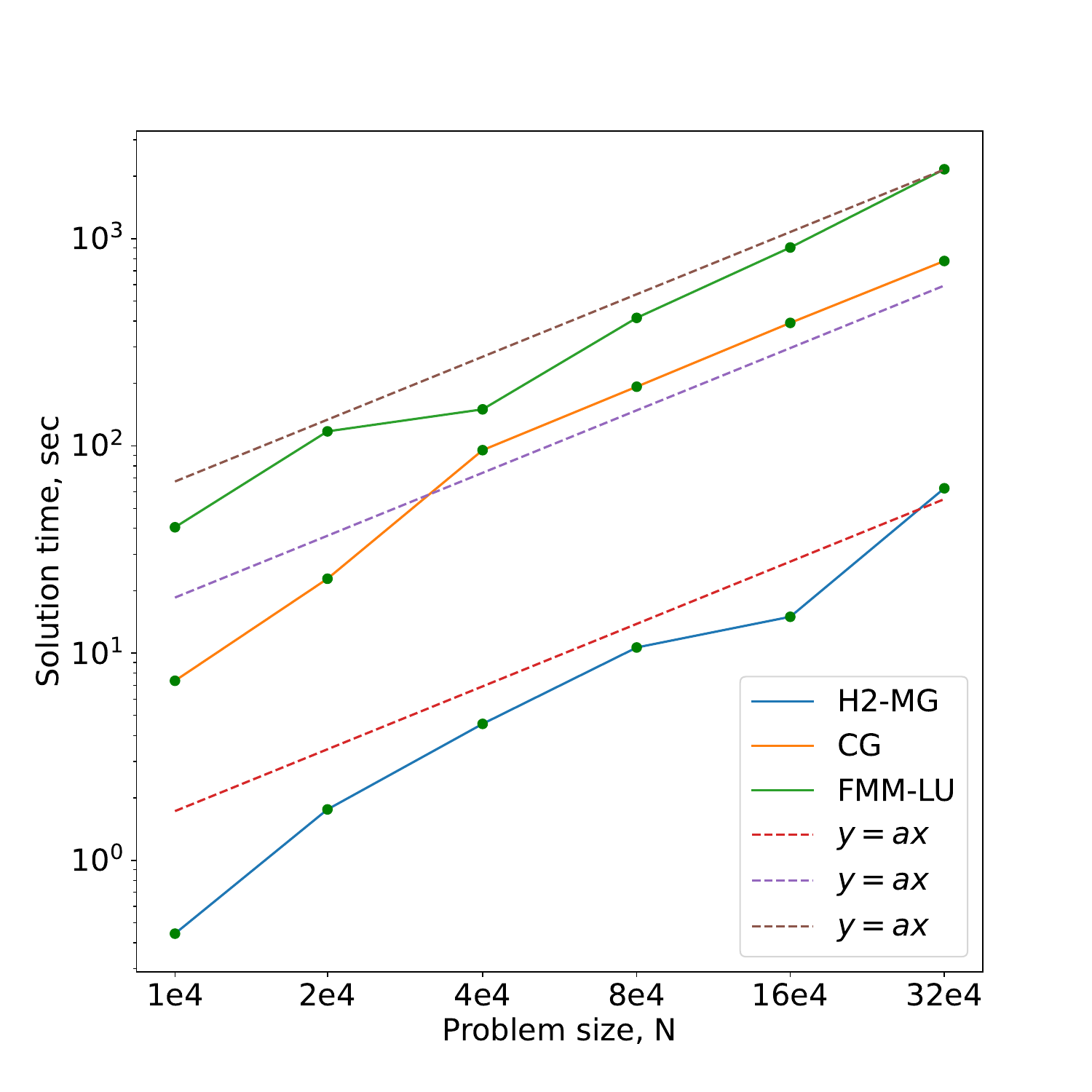}}
     \subcaptionbox{ $c = 10^{-5}$ \label{fig:gau_st_s01_5}}{
        \includegraphics[width=0.48\textwidth]{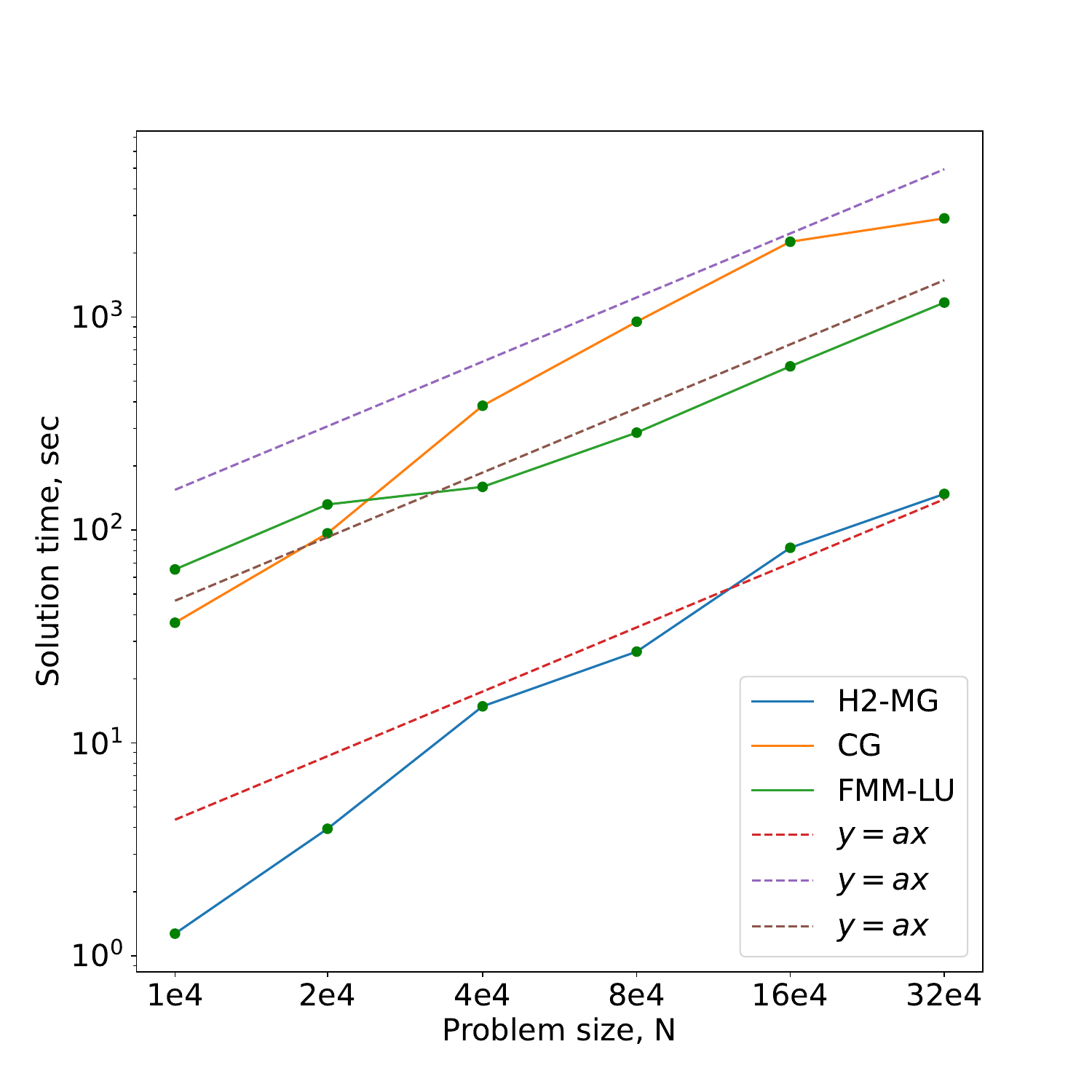}}
        \caption{Overall solution time of methods \HTMG and CG,  $\sigma = 0.1$ }
        \label{fig:gau_nvc_s01_38}
\end{figure}

\begin{figure}[h]
    \subcaptionbox{ $c = 10^{-3}$ \label{fig:gau_st_s001_3}}{
        \includegraphics[width=0.49\textwidth]{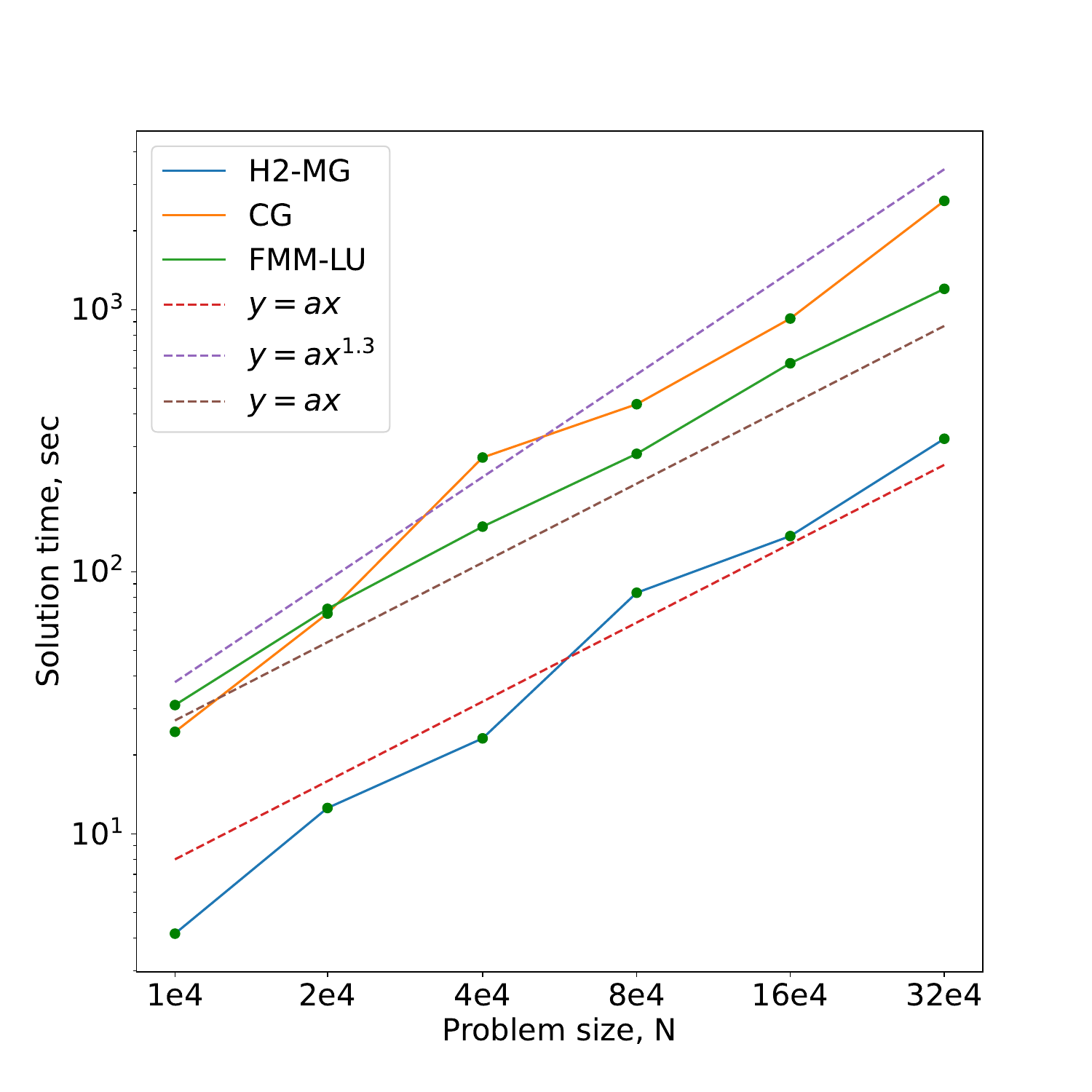}}
    \subcaptionbox{ $c = 10^{-5}$ \label{fig:gau_st_s001_5}}{
        \includegraphics[width=0.49\textwidth]{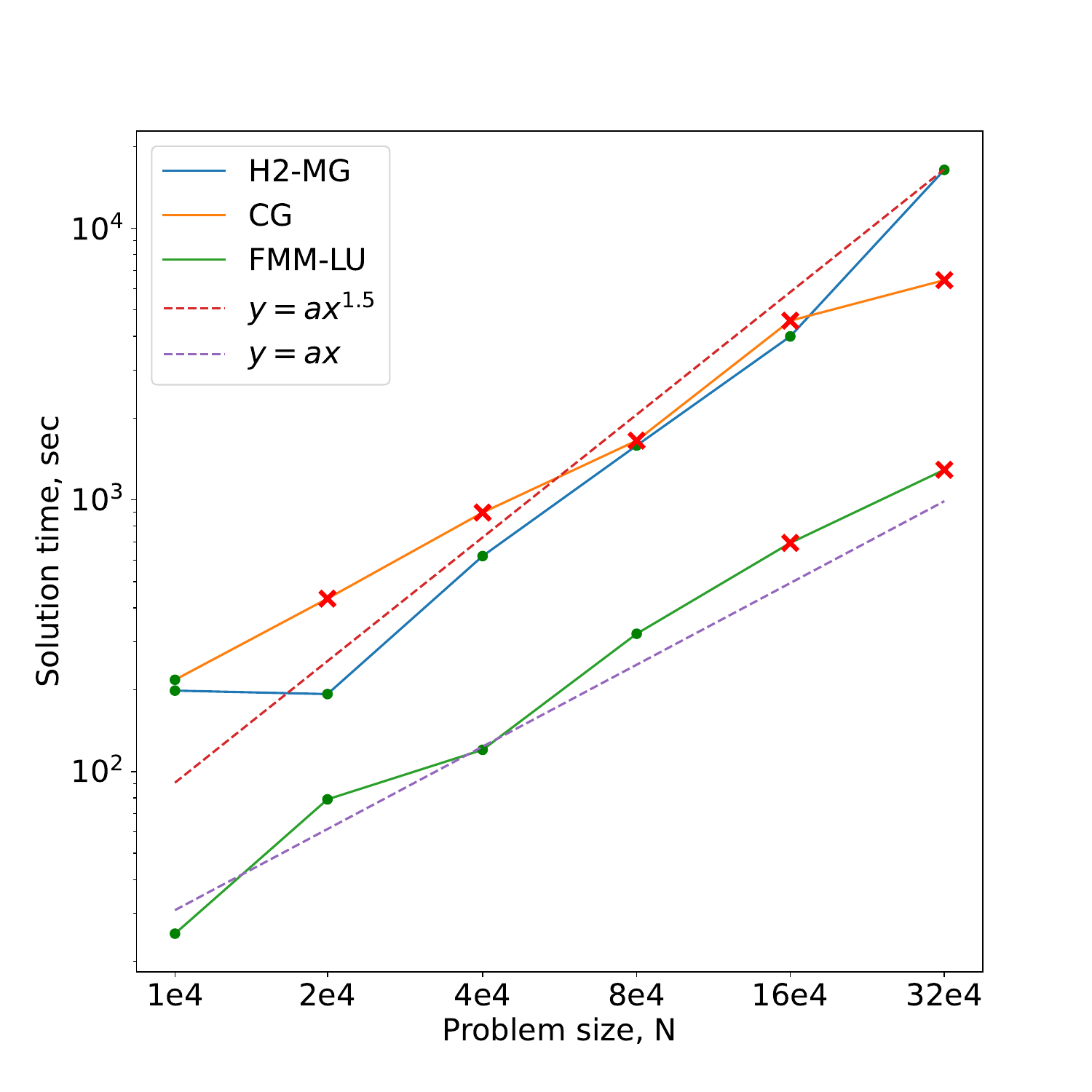}}
        \caption{Overall solution time of methods \HTMG and CG,  $\sigma = 0.01$}
        \label{fig:gau_st_s001}
\end{figure}

Figures~\ref{fig:gau_nvc_s01_38} and~\ref{fig:gau_st_s001}  show the comparison of the overall solution time for the \HTMG and CG methods, compared against the H2-direct solver FMM-LU~\cite{sush-fmmlu-2023}. The goal of the comparison is to explore the efficiency of the proposed method by benchmarking it against an alternative efficient \HT solver. A red cross indicates that either the iterative method failed to converge within 5000 iterations, or the direct solver failed to solve the system with the required accuracy.

For $\sigma=0.01$ and $c=10^{-5}$, the system appears to have an extremely large condition number, leading to superlinear scaling or failure to converge within 5000 iterations of the iterative solvers. For the direct solver, this condition results in failure due to computational errors, as LU methods without pivoting struggle to handle systems with extremely ill-conditioned matrices. For three other cases, both \HTMG and FMM-LU solve the system with linear scaling. Since both \HTMG and FMM-LU scale linearly in this example, the main competition lies in the constant factor, where the iterative \HTMG naturally outperforms direct solver FMM-LU. \HTMG also outperforms CG in constant.




\subsection{Exponential kernel}
\label{seq:num_exp}
Consider the linear system~\eqref{eq:axb_num} where $b\in\R{N}$ is a right-hand side vector, $x\in\R{N}$ is an unknown vector, and $A\in\R{N\times N}$ is a kernel matrix with the linear exponential decay kernel. 
Consider a uniform tensor grid on a unit square $P\subset \R{2}$: $p_i\in P$, $i\in 1\dots N$, 
where $N$ is the number of points. The kernel matrix $A$ is given by the formula:
\begin{equation}
\label{eq:a_exp}
   a_{ij} = \begin{cases}
\exp({-\frac{|p_i-p_j|}{\sigma}}),          & \quad \text{if } i\neq j\\
1 + c,                                        & \quad \text{if } i=j
\end{cases}, 
\end{equation}
\noindent
where $\sigma\in\R{}$ is the dispersion parameter of the matrix, $c\in \R{}$ is a constant. Matrix $A$ is approximated to the \HT format with accuracy $\epsilon=10^{-9}$, number of levels is chosen adaptively.  

We keep all the assumptions made for the problem in Section~\ref{seq:num_gau} for this example.

\subsubsection{Exponential kernel matrix, asymptotics analysis}

\begin{figure}[htb]
     \centering
     \subcaptionbox{ $c = 10^{-3}$ \label{fig:exp_nvc_s01}}{
        \includegraphics[width=0.48\textwidth]{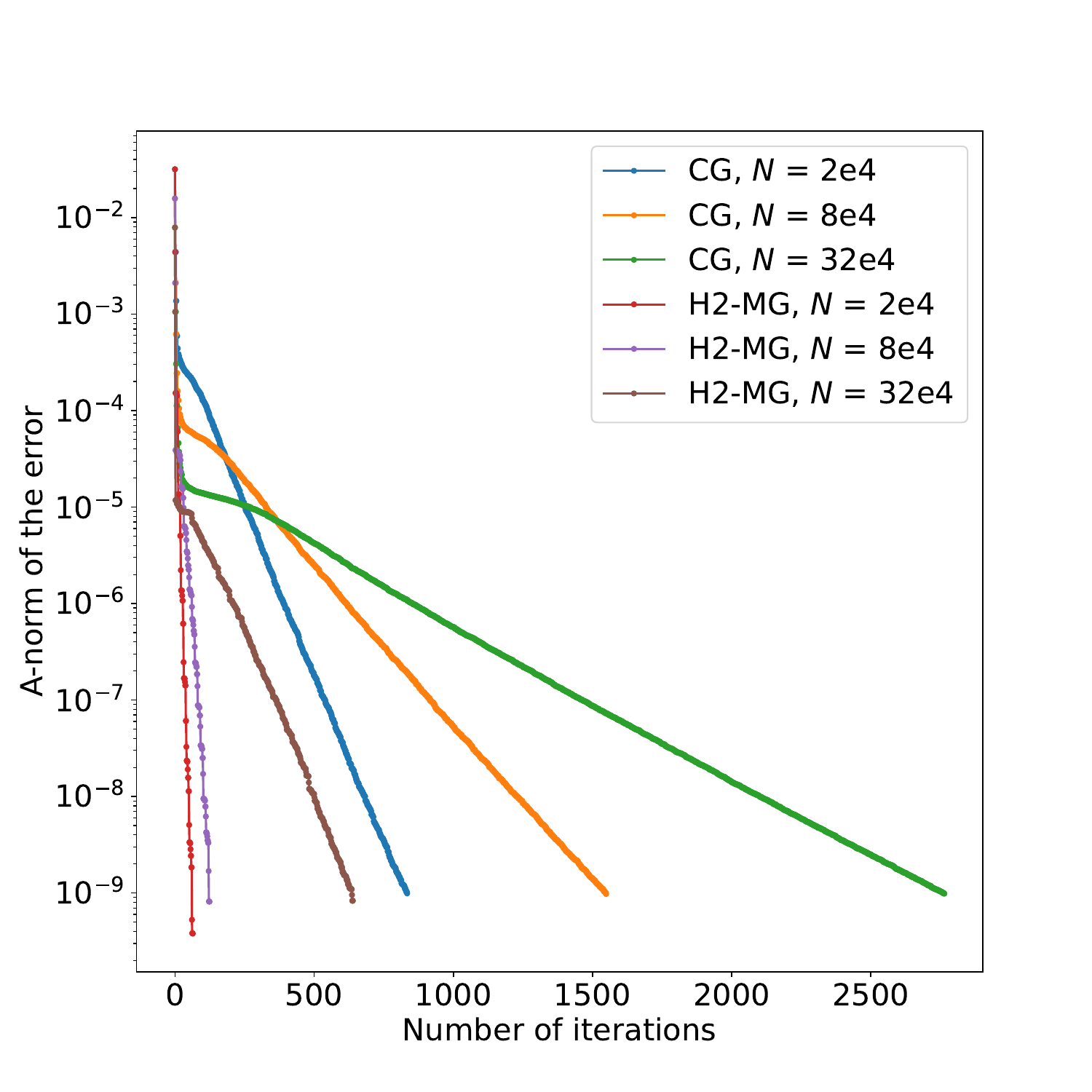}}
     \subcaptionbox{$c = 10^{-5}$ \label{fig:exp_nvc_s01_2}}{
        \includegraphics[width=0.48\textwidth]{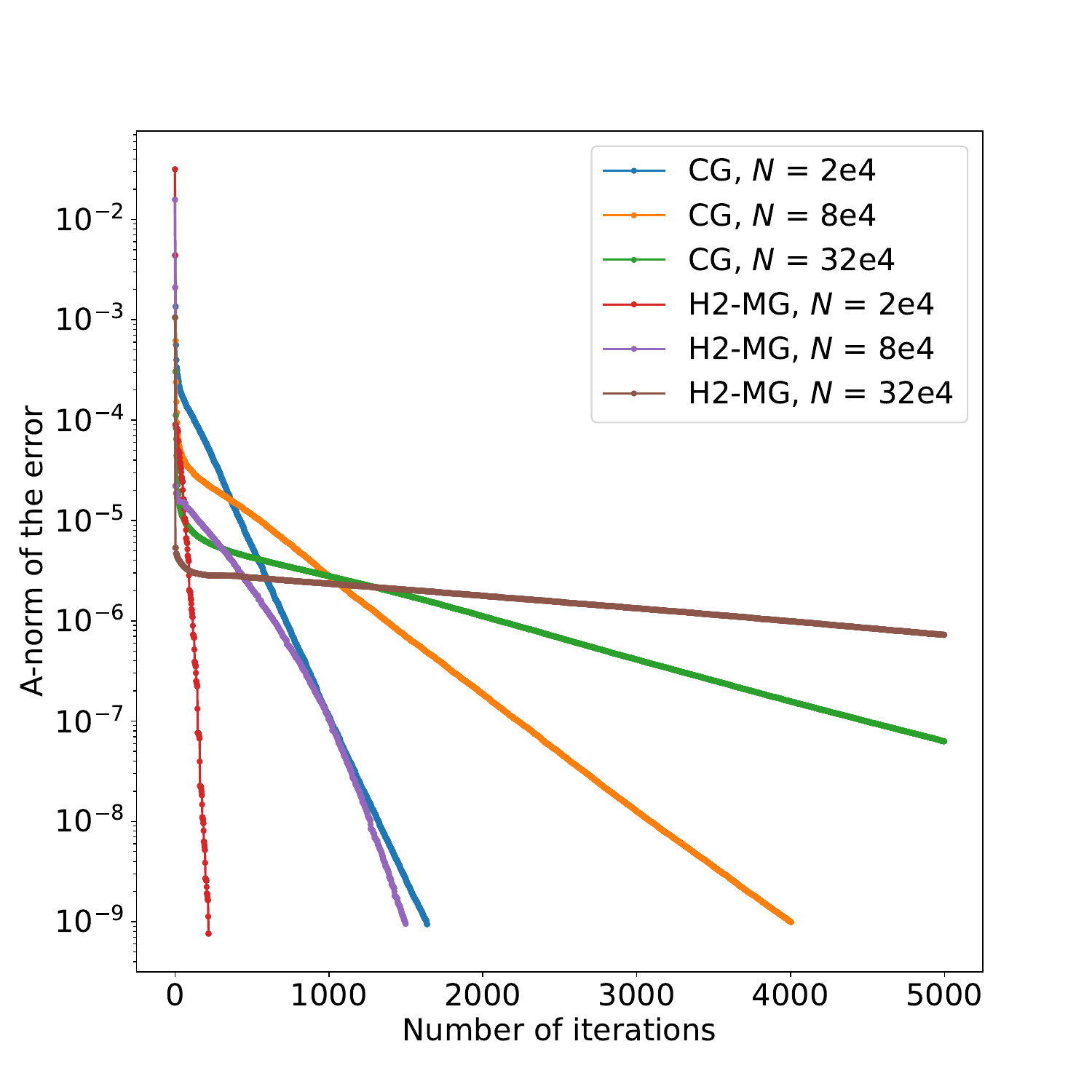}}
        \caption{Convergence evolution of \HTMG and CG as problem size increases,  $\sigma = 0.1$}
        \label{fig:exp_nvc_s01_16}
\end{figure}

\begin{figure}[htb]
    \subcaptionbox{ $c = 10^{-3}$ \label{fig:exp_nvc_s001}}{
        \includegraphics[width=0.48\textwidth]{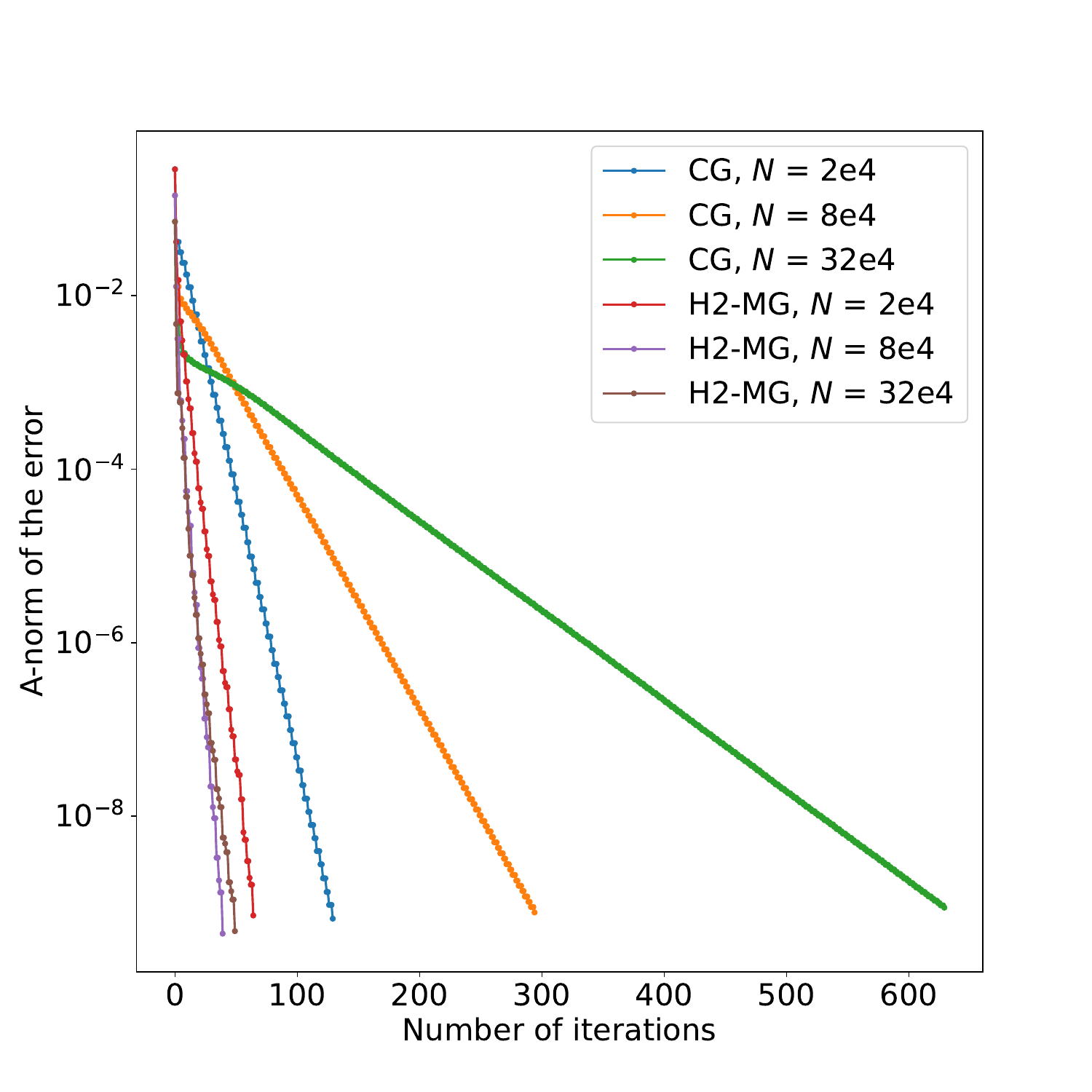}}
    \subcaptionbox{ $c = 10^{-5}$ \label{fig:exp_nvc_s01_34}}{
        \includegraphics[width=0.48\textwidth]{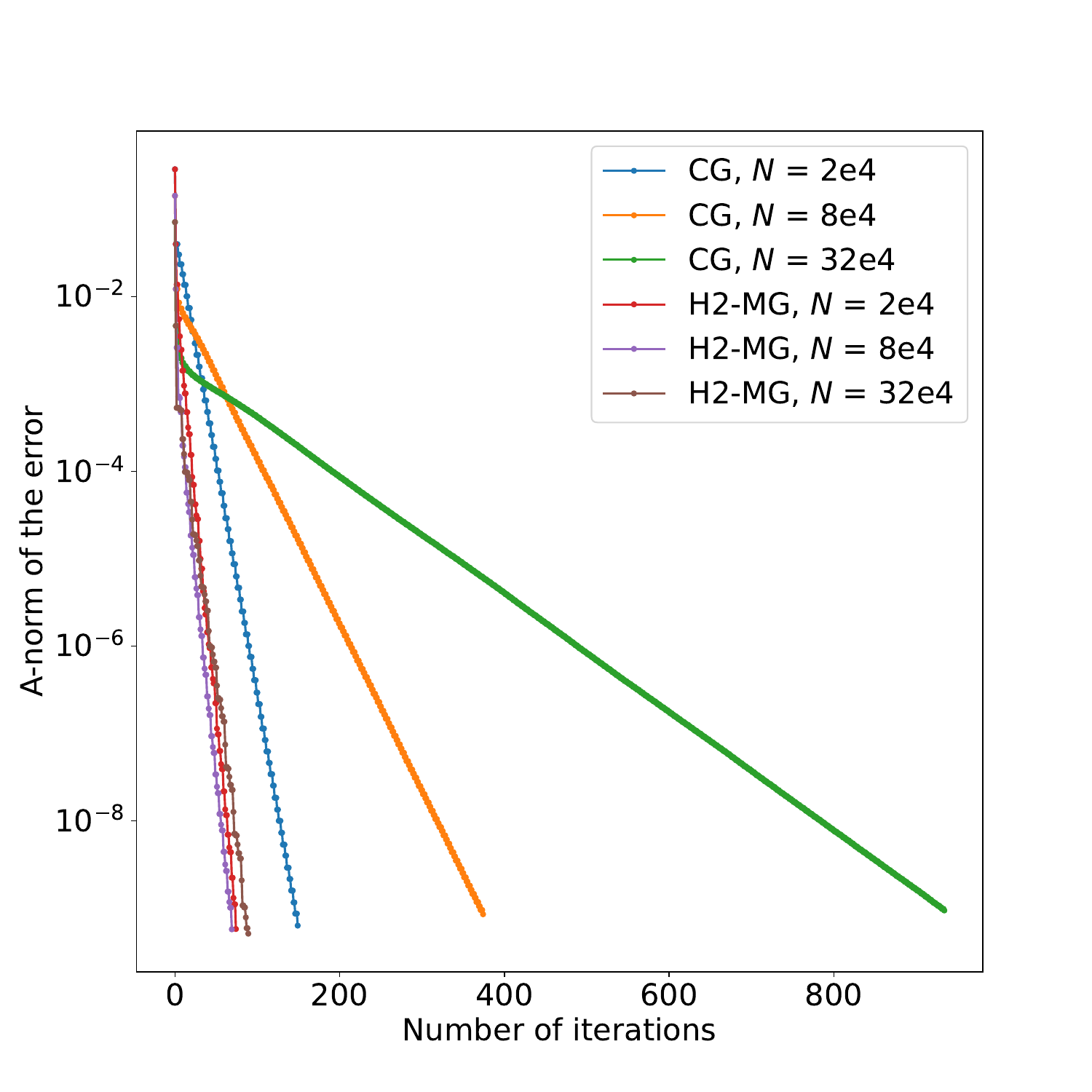}}
        \caption{Convergence evolution of \HTMG and CG as problem size increases,  $\sigma = 0.01$}
        \label{fig:exp_nvc_s001_17}
\end{figure}

In this subsection we consider the asymptotic behavior of the system~\eqref{eq:axb_num} with matrix $A$ given by~\eqref{eq:a_exp} for various combinations of parameters $c$ and $\sigma$.
Figures~\ref{fig:exp_nvc_s01_16}~and~\ref{fig:exp_nvc_s001_17} show the convergence of the \HTMG algorithm compared to CG across various problem sizes. 
The method tolerance is  $\varepsilon = 10^{-9}$.
We use parameters $n_f = 1$, $n_c = 40$.

We can see that the strong regularization parameter $c=10^{-3}$ leads to the fast convergence of the \HTMG method for both matrices with $\sigma = 0.1$ and $\sigma = 0.01$, while the weak regularization parameter $c=10^{-5}$ leads to the divergence of both algorithms for the larger problem sizes. 

The numbers of V-cycles required for the convergence to the fixed accuracy are presented in Table~\ref{tab:exp_d2_vc}.

\begin{table}[htb]
\centering
\begin{tabular}{c||c|c|c|c|c|c}
    \multirow{2}{*}{Matrix Parameters} & \multicolumn{6}{c}{Problem Size} \\ \cline{2-7}
    & 1e4 &  2e4   &  4e4 & 8e4 &  16e4 &  32e4  \\ \hline
   $\sigma = 0.1$, $c = 10^{-3}$  & 14 &  13   &  18  & 25   &  54   &  128     \\
   $\sigma = 0.1$, $c = 10^{-5}$  & 31 &  44    &  127  & 300  &  -   &  -     \\
   $\sigma = 0.01$, $c = 10^{-3}$  & 19 &  13  &  25  & 8   &  19  &  10    \\
   $\sigma = 0.01$, $c = 10^{-5}$  & 22 &  15   &  30 &  14   &  30  &  18    \\
\end{tabular}
\caption{Number of V-cycles to solve the system for different problem sizes}
\label{tab:exp_d2_vc}
\end{table}

As in the previous example, the number of V-cycles required for the convergence to the fixed accuracy does not grow significantly.

\begin{figure}[H]
     \centering
     \subcaptionbox{ $c = 10^{-3}$ \label{fig:exp_st_s01_1}}{
        \includegraphics[width=0.48\textwidth]{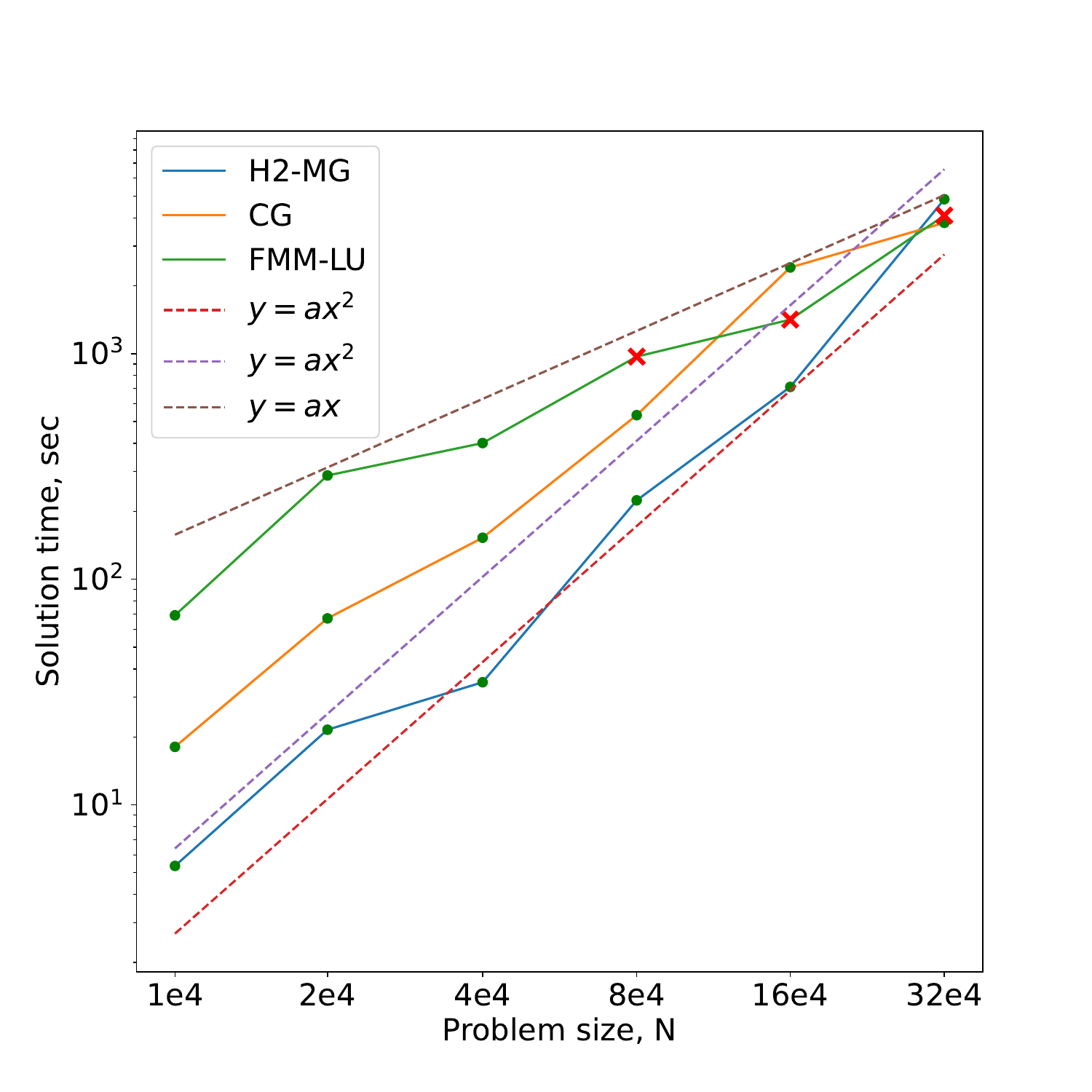}}
     \subcaptionbox{ $c = 10^{-5}$ \label{fig:exp_st_s01_2}}{
        \includegraphics[width=0.48\textwidth]{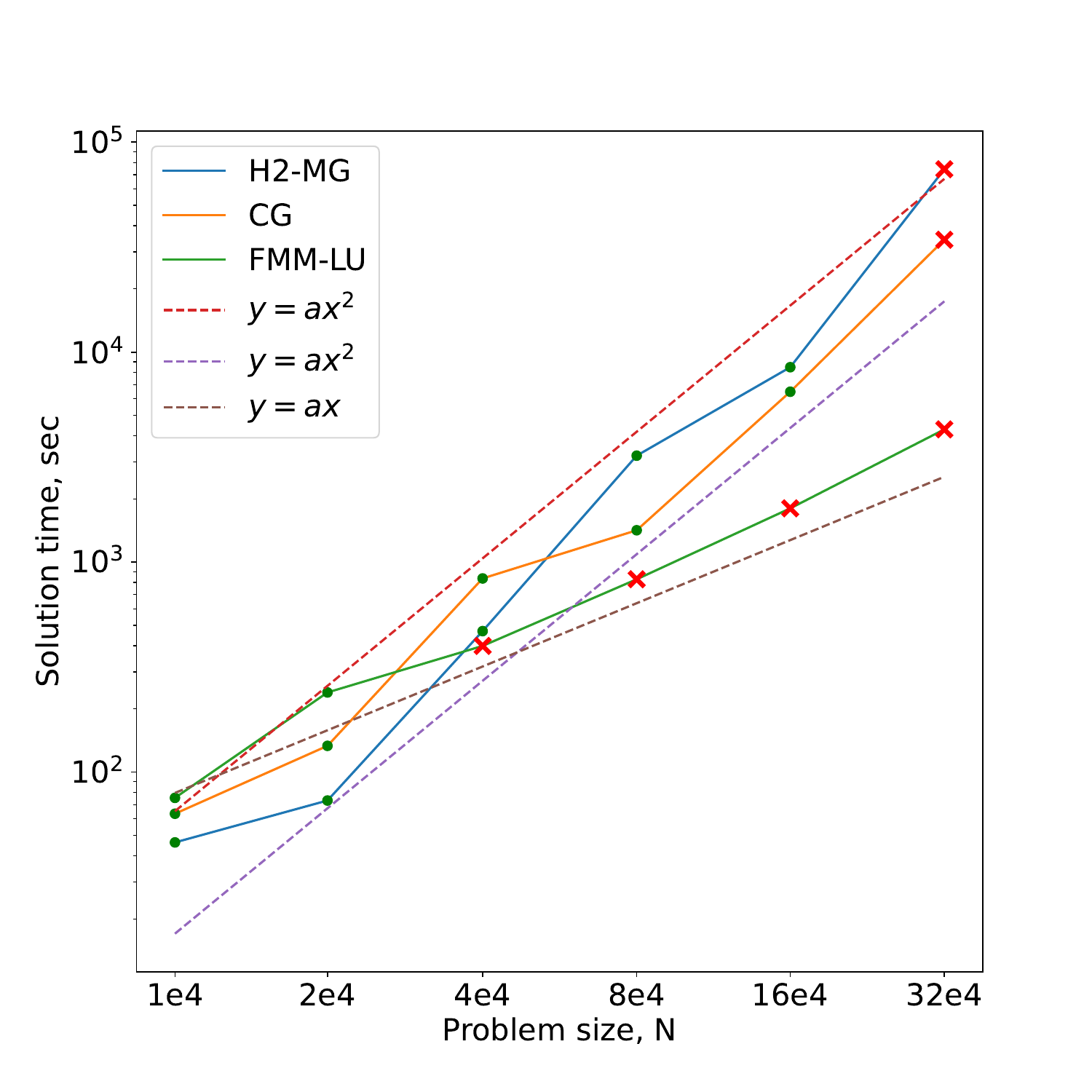}}
        \caption{Overall solution time of methods \HTMG and CG,  $\sigma = 0.1$ }
        \label{fig:exp_nvc_s01_3}
\end{figure}

\begin{figure}[H]
    \subcaptionbox{ $c = 10^{-3}$ \label{fig:exp_st_s01_3}}{
        \includegraphics[width=0.49\textwidth]{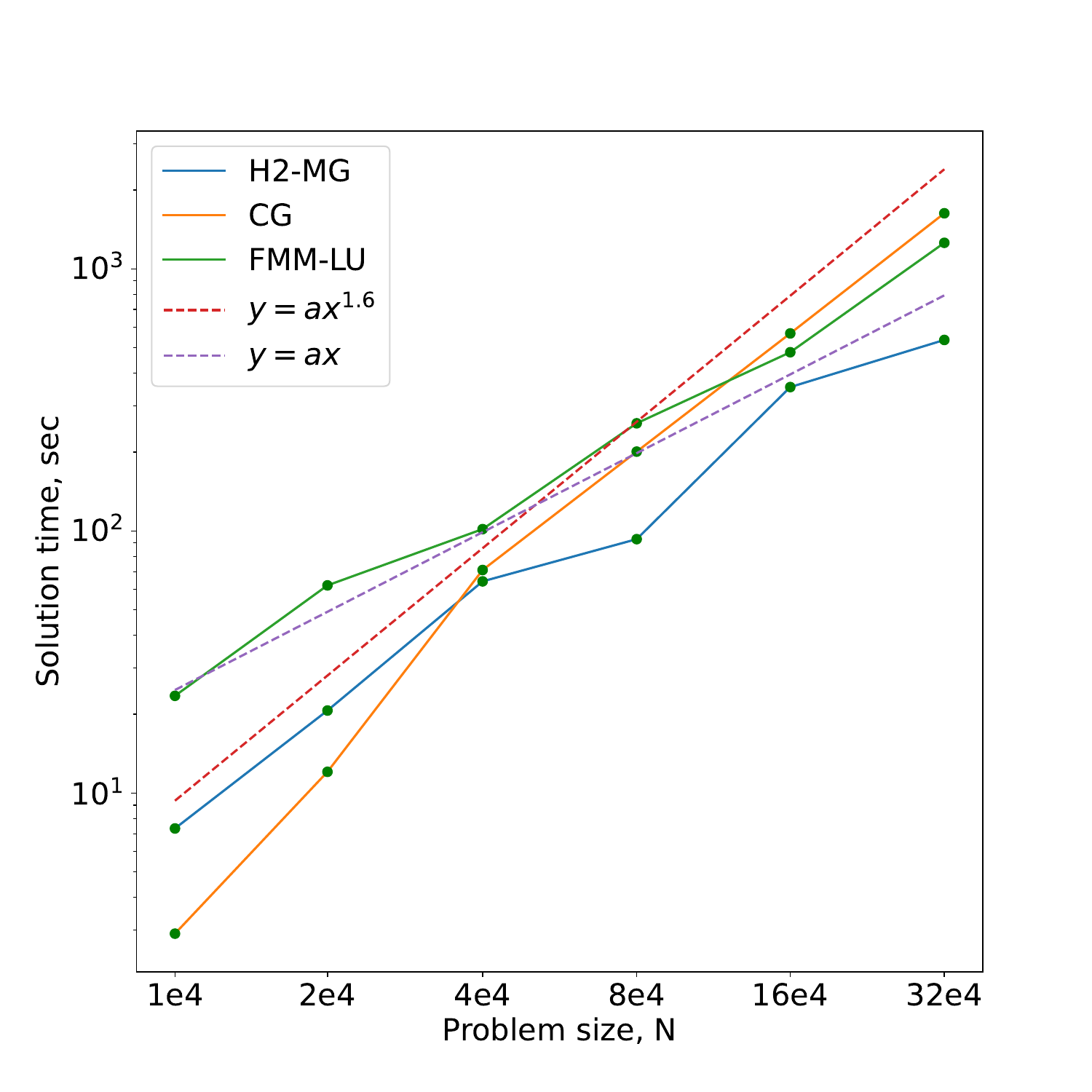}}
    \subcaptionbox{ $c = 10^{-5}$ \label{fig:exp_st_s01_4}}{
        \includegraphics[width=0.49\textwidth]{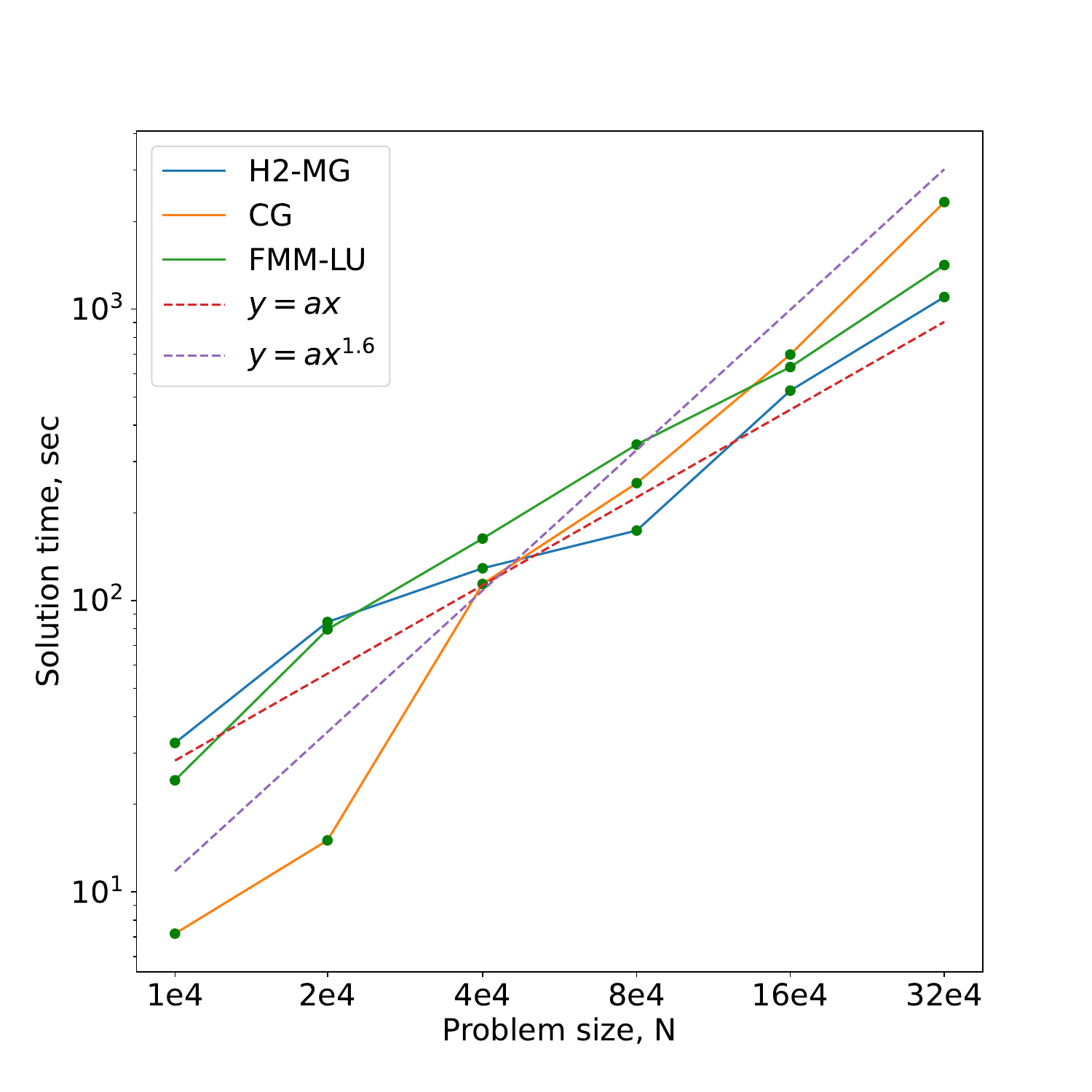}}
        \caption{Overall solution time of methods \HTMG and CG,  $\sigma = 0.01$}
        \label{fig:exp_st_s001}
\end{figure}

Figures~\ref{fig:exp_nvc_s01_3} and~\ref{fig:exp_st_s001} show the comparison of the overall solution time for the \HTMG and CG methods, compared against the H2-direct solver FMM-LU. A red cross {\color{black} in Figure~\ref{fig:exp_nvc_s01_3}} indicates that either the iterative method failed to converge within 5000 iterations, or the direct solver failed to solve the system with the required accuracy.

For $\sigma=0.1$, the system appears to have an extremely large condition number, leading to quadratic scaling or failure to converge within 5000 iterations of the iterative solvers. For the direct solver, this condition results in failure due to computational errors, as LU methods without pivoting struggle to handle systems with extremely ill-conditioned matrices. For $\sigma=0.01$, both \HTMG and FMM-LU solve the system with linear scaling, while CG exhibits a higher scaling power of $y = x^{1.5}$. Since both \HTMG and FMM-LU scale linearly in this example, the main competition lies in the constant factor, where the iterative \HTMG naturally outperforms.

{\color{black}
\subsection{Boundary element method for a 3D electrostatic problem}
For the third example, we consider Laplace’s equation in its integral form, solved using the Boundary Element Method (BEM) on a complex 3D surface. This formulation is widely used in electrostatics, capacitance computations, molecular solvation models (such as the Poisson–Boltzmann equation), and modeling interactions between charged surfaces.

Specifically, we consider the \textit{Laplace single-layer potential}, which describes the potential $u(\mathbf{x})$ at a point $\mathbf{x} \in \Gamma \subset \mathbb{R}^3$, where $\Gamma$ is a given surface. The potential $u(\mathbf{x})$ arises from charges (or equivalent sources) $\sigma(\mathbf{y})$ on the same surface $\Gamma$. The Laplace single-layer potential formula is:
$$
u(\mathbf{x}) = \int_\Gamma \frac{1}{4\pi \|\mathbf{x} - \mathbf{y}\|} \, \sigma(\mathbf{y}) \, dS_{\mathbf{y}},
$$
where $\|\mathbf{x} - \mathbf{y}\|$ is the Euclidean distance between source and observation points, and $dS_{\mathbf{y}}$ is the surface area element.

To solve this numerically, we discretize the surface $\Gamma$ into $N$ triangular elements with associated centroids $\{ \mathbf{y}_j \}$, quadrature weights $w_j$, and evaluation points $\{ \mathbf{x}_i \}$. 

To account for the singularity of the kernel when $i = j$, we modify the diagonal entries using a geometrically motivated regularization. We define $R_i$ as the average distance from the centroid of triangle $i$ to its three vertices. The resulting matrix $A \in \mathbb{R}^{N \times N}$ represents the discretized integral operator, with entries defined as:

\begin{equation}
    A_{ij} =
\begin{cases}
\displaystyle \frac{w_j}{4\pi \|\mathbf{x}_i - \mathbf{y}_j\|}, & \text{if } i \neq j, \\[10pt]
\displaystyle \frac{w_i}{4\pi R_i}, & \text{if } i = j,
\end{cases}
\end{equation}
where
$$
R_i = \frac{1}{3} \sum_{k=1}^3 \| \mathbf{y}_i - \mathbf{v}_{ik} \|,
$$
and $\{ \mathbf{v}_{i1}, \mathbf{v}_{i2}, \mathbf{v}_{i3} \}$ are the vertices of triangle $i$.

For our numerical experiment, we chose a complex 3D geometry: a curved torus shown in Figure~\ref{fig:wtorus_sigma}. 
To construct the right-hand side of the system, we consider a point source located at  $\mathbf{x}_0 \in \mathbb{R}^3$ outside the surface $\Gamma$. The potential generated by this point source at a location $\mathbf{y} \in \Gamma$ is given by the free-space Green’s function:
$$f(\mathbf{y}) = \frac{1}{4\pi \| \mathbf{y} - \mathbf{x}_0 \|}.$$
We evaluate this expression at each centroid $\mathbf{y}_j$ of the surface elements to obtain the right-hand side vector $\mathbf{f} \in \mathbb{R}^N$:
$$f_j = \frac{1}{4\pi \| \mathbf{y}_j - \mathbf{x}_0 \|}, \quad j = 1, \dots, N.$$
 We compare the CG method with \HTMG method on this problem. The comparison is performed for several mesh resolutions on the curved torus geometry. We set the number of fine grid iterations to $n_f = 1$ and the number of coarse grid iterations to $n_c = 40$. Figures~\ref{fig:bem_iter} and~\ref{fig:bem_time} show the convergence behavior of both methods in terms of the number of iterations and total computational time, respectively. In this example, we plot the residual, since the exact solution is unknown and we cannot compute the $A$-norm of the error, as we did in the previous examples.
\begin{figure}[h]
     \centering
     \subcaptionbox{\color{black}Residual per fine-level iteration \label{fig:bem_iter}}{
        \includegraphics[width=0.48\textwidth]{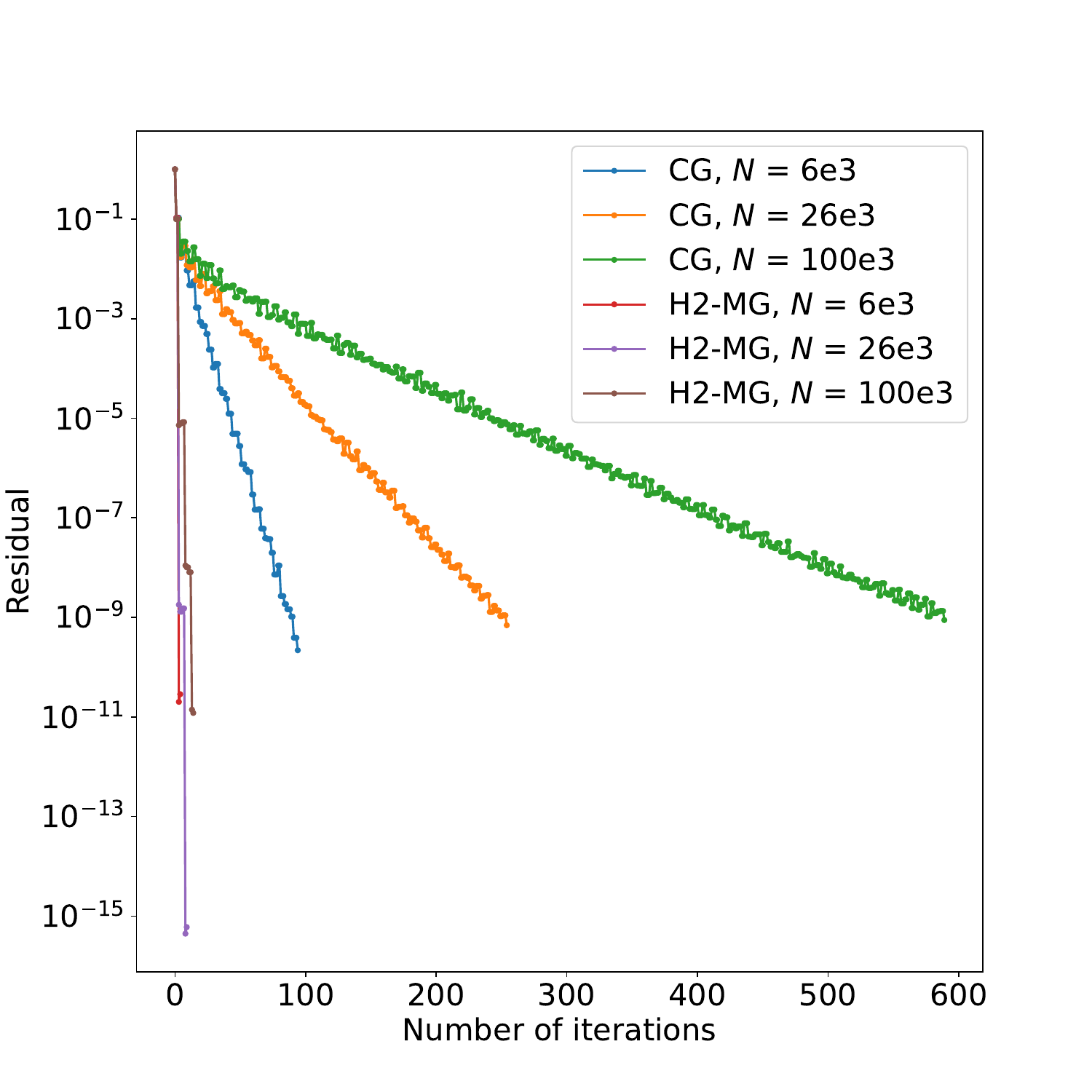}}
     \subcaptionbox{\color{black}Residual per time \label{fig:bem_time}}{
        \includegraphics[width=0.48\textwidth]{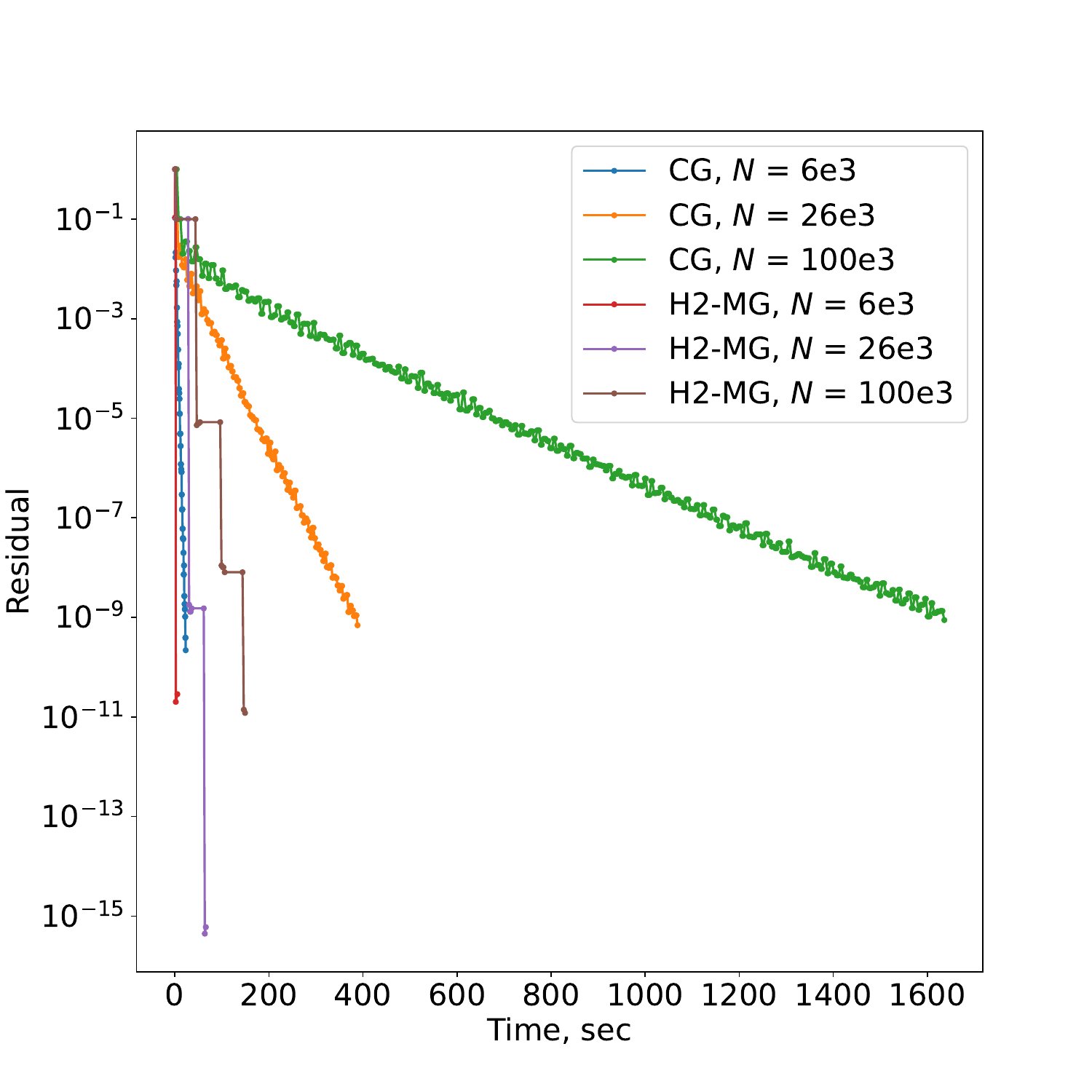}}
        \caption{\color{black}Convergence evolution of \HTMG and CG as problem size increases}
        \label{fig:bem_conv}
\end{figure}
We can see from the Figure~\ref{fig:bem_iter} that the number of iterations required to solve the system using \HTMG does not increase significantly with problem size, while for CG it grows substantially. To confirm this effect quantitatively, Table~\ref{tab:bem_vc} reports the number of V-cycles for \HTMG and the number of iterations for CG. Note that these numbers should not be compared directly, as one V-cycle includes several inner CG iterations. Instead, we focus on the growth trend.
\begin{table}[h]
    \color{black}
    \centering
    \begin{tabular}{c||c|c|c|c|c}
    \multirow{2}{*}{Method} & \multicolumn{5}{c}{Problem Size} \\ \cline{2-6}
     & 6e3 & 12e3 & 26e3 & 50e3 & 100e3 \\ \hline
    \HTMG & 1 & 2 & 2 & 2 & 3 \\ \hline
    CG    & 95 & 155 & 255 & 380 & 590 \\ 
\end{tabular}
    \caption{\color{black}Comparison of V-Cycles of \HTMG and iterations of CG across problem sizes}
    \label{tab:bem_vc}
\end{table}

Figure~\ref{fig:sol_time_bem} shows the total solution time for the system using CG and \HTMG.
We observe that \HTMG outperforms CG in timing not only in terms of the constant factor but also in asymptotic scaling with problem size.
Figure~\ref{fig:wtorus_sigma} presents the computed surface charge density $\sigma$ on the curved torus $\Gamma$.
This example demonstrates that the \HTMG algorithm can be effectively applied to BEM problems.

\begin{figure}[H]  
\centering
        \includegraphics[width=0.5\textwidth]{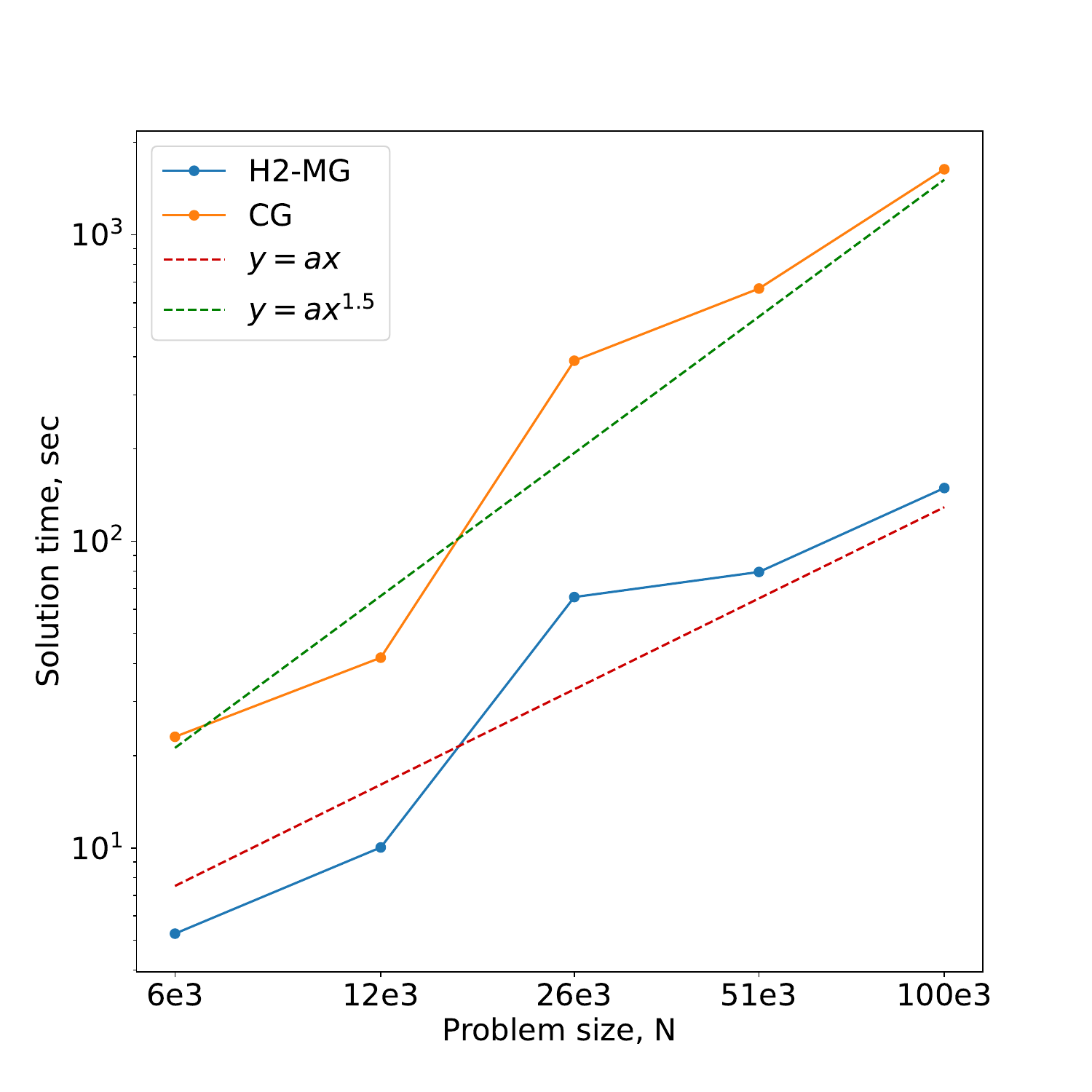}
        \caption{\color{black} Total solution time for CG and \HTMG across mesh sizes}
        \label{fig:sol_time_bem}
\end{figure}

\begin{figure}[H]  
\centering
        \includegraphics[width=0.7\textwidth]{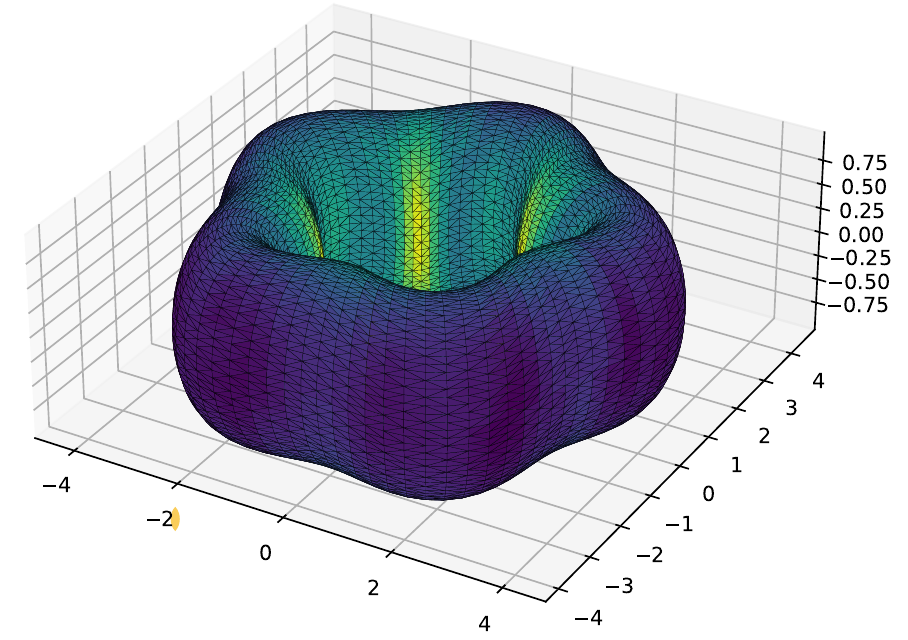}
        \caption{\color{black}Surface charge $\sigma$ distribution on the triangulated surface $\Gamma$ of a curved torus}
        \label{fig:wtorus_sigma}
\end{figure}
}

\section{Conclusion}
 The \HTMG algorithm offers an advance in solving large, dense kernel matrices efficiently by iterative means. By combining the rapid convergence of the multigrid method with the time and memory efficiencies of \HT matrix approximations, this algorithm not only fills a gap in the existing suite of \HT solvers but also expands the toolkit available for tackling complex computational problems. The demonstrated linear complexity and practical effectiveness of $\mathcal{H}^{2}$-MG, verified through numerical examples, underscores its potential for applications burdened by large, dense, kernel matrices. Future work will aim to expand its applicability beyond symmetric positive definite matrices, high-performance implementation, and integration into other computational frameworks.
\bibliographystyle{siamplain}
\bibliography{lib}
\end{document}